\documentclass[11pt]{article}
\usepackage{latexsym}
\usepackage{amsmath}
\usepackage{cite}
\usepackage{amsthm,color,bm}
\usepackage{amssymb}
\usepackage{mathrsfs}
\usepackage{indentfirst}
\usepackage{times,cases}
\usepackage{dsfont}

\usepackage{lscape}
\newtheorem{theorem}{\noindent Theorem}[section]
\newtheorem{proposition}[theorem]{\noindent Proposition}
\newtheorem{definition}[theorem]{\noindent Definition}
\newtheorem{lemma}[theorem]{\noindent Lemma}
\newtheorem{remark}[theorem]{\noindent Remark}
\newtheorem{corollary}[theorem]{\noindent Corollary}
\numberwithin{figure}{section}
\numberwithin{equation}{section}

\renewcommand{\theequation}{\thesection.\arabic{equation}}

\DeclareMathOperator*{\esssup}{ess\,sup}
\DeclareMathOperator*{\LIM}{LIM}

\newcommand{\cA}{{\mathcal A}}
\newcommand{\cB}{{\mathcal B}}
\newcommand{\cC}{{\mathcal C}}
\newcommand{\cD}{{\mathcal D}}
\newcommand{\cE}{{\mathcal E}}
\newcommand{\cF}{{\mathcal F}}
\newcommand{\cG}{{\mathcal G}}

\newcommand{\cI}{{\mathcal I}}

\newcommand{\cO}{{\mathcal O}}

\newcommand{\cR}{{\mathcal R}}

\newcommand{\cT}{{\mathcal T}}
\newcommand{\cU}{{\mathcal U}}

\newcommand{\cZ}{{\mathcal Z}}

\newcommand{\sA}{{\mathscr A}}

\newcommand{\sM}{{\mathscr M}}
\newcommand{\sN}{{\mathscr N}}

\newcommand{\fA}{{\mathfrak A}}
\newcommand{\fU}{{\mathfrak U}}

\def\T{\mathbb{T}}

\def\Z{\mathbb{Z}}
\def\R{\mathbb{R}}

\def\bP{\mathbb{P}}

\def\N{\mathbb{N}}
\def\1{\mathds{1}}
\newcommand{\di}{\mathrm{d}}
\newcommand{\me}{\mathrm{e}}

\def\disp{\displaystyle}

\def\bc{\begin{center}}
	\def\ec{\end{center}}
\def\be{\begin{equation}}
	\def\ee{\end{equation}}
\def\bea{\begin{eqnarray}}
	\def\eea{\end{eqnarray}}
\def\ba{\begin{array}}
	\def\ea{\end{array}}
\def\benu{\begin{enumerate}}
	\def\eenu{\end{enumerate}}
\def\bt{\begin{theorem}}
	\def\et{\end{theorem}}
\def\bl{\begin{lemma}}
	\def\el{\end{lemma}}
\def\bco{\begin{corollary}}
	\def\eco{\end{corollary}}
\def\bn{\begin{numcases}}
	\def\en{\end{numcases}}
\def\br{\begin{remark}}
	\def\er{\end{remark}}
\def\bd{\begin{definition}}
	\def\ed{\end{definition}}
\def\bp{\begin{proposition}}
	\def\ep{\end{proposition}}
\def\bo{\begin{proof}}
	\def\eo{\end{proof}}
\def\bx{\begin{example}}
	\def\ex{\end{example}}
\def\bal{\begin{align}}
	\def\eal{\end{align}}

\def\pa{\partial}
\def\al{\alpha}
 \def\de{\delta}

\def\lam{\lambda} \def\Lam{\Lambda}

\def\ve{\varepsilon}
\def\sig{\sigma}
\def\vsig{\varsigma}
\def\vp{\varphi}

\def\w{\omega}\def\W{\Omega}

\def\~{\widetilde}

\def\ol{\overline}

\def\Cap{\bigcap}
\def\Cup{\bigcup}

\def\ra{\rightarrow}
\def\Ra{\Rightarrow}

\def\8{\infty}
\def\X{\times}

\def\mb{\mbox}
\def\di{{\rm d}}
\def\me{{\rm e}}

\def\suo{\!\!\!}

\def\Hs{\hspace{0.8cm}}
\def\hs{\hspace{0.4cm}}
\def\Vs{\vskip10pt}
\def\vs{\vskip5pt}

\def\({\left(}
\def\){\right)}

\begin{document}

	
	\begin{center}
		{\large \bf Convergence of bi-spatial pullback random attractors and
stochastic Liouville type\\equations for nonautonomous stochastic $p$-Laplacian lattice system}
		\vspace{0.5cm}\\
		{Jintao Wang*,\hs Qinghai Peng,\hs Chunqiu Li}\\\vspace{0.3cm}
		
		{\small Department of Mathematics, Wenzhou University, Wenzhou 325035, China}
	\end{center}

	
	\renewcommand{\theequation}{\arabic{section}.\arabic{equation}}
	\numberwithin{equation}{section}


\begin{abstract}
We consider convergence properties of the long-term behaviors with respect to the coefficient of
the stochastic term for a nonautonomous stochastic $p$-Laplacian lattice equation with multiplicative noise.
First, the upper semi-continuity of pullback random $(\ell^2,\ell^q)$-attractor is
proved for each $q\in[1,+\8)$.
Then, a convergence result of the time-dependent invariant sample Borel probability measures
is obtained in $\ell^2$.
Next, we show that the invariant sample measures satisfy a stochastic Liouville type equation
and a termwise convergence of the stochastic Liouville type equations is verified.
Furthermore, each family of the invariant sample measures is turned out to be a sample statistical solution,
which hence also fulfills a convergence consequence.
\vs
		
\vs
		
		
\end{abstract}
	
	\vspace{-1 cm}
	
\footnote[0]{\hspace*{-7.4mm}
$^{*}$ Corresponding author.\\
E-mail address: wangjt@wzu.edu.cn (J.T. Wang); 22451025013@stu.wzu.edu.cn (Q.H. Peng); lichunqiu@wzu.edu.cn (C.Q.
Li).}

\section{Introduction}

This article investigates convergence properties of the long-term dynamical behaviors
(as $\al$ changes) of the following nonautonomous stochastic $p$-Laplacian
lattice equation with multiplicative noise,
\begin{align}\label{1.1}
\di u_i(t)=&\nu(t)\left[|u_{i+1}-u_i|^{p-2}(u_{i+1}-u_i)
-|u_i-u_{i-1}|^{p-2}(u_i-u_{i-1})\right]\di t\notag\\
&+\(f_i(t,u_i)-\lam u_i\)\di t+\alpha u_i\di W(t),\hs t>\tau
\end{align}
with the initial datum
\be\label{1.2}u_i(\tau)=u_{\tau i}\hs\mb{and}\hs i\in\Z,\ee
where the coefficient $\nu(t)\in\cC(\R)$ is assumed to satisfy
$0\leqslant\nu(t)\leqslant\nu_0<\8$, $p\geqslant2$, $f$ is the nonautonomous nonlinear forcing,
$W$ is a two-sided real-valued Wiener process on the probability space to be determined later
and $u_{\tau}:=(u_{\tau i})_{i\in\Z}\in\ell^2$.
The lattice equation \eqref{1.1} can be regarded as the spatial discrete form
on infinite lattices of 1-dimensional $p$-Laplacian equation (see \cite{GL17})
\be\label{1.3}
\frac{\pa u}{\pa t}=\nu(t)(|u_x|^{p-2}u_x)_x+f(t,u)-\lam u+\al u\di W(t),\hs x\in\R.\ee

The $p$-Laplacian equation originated from the study of a wide variety of physical phenomena
(\cite{GK16,WZ23}) involving non-Newtonian fluids, porous media and nonlinear elasticity etc.
There have been plenty of works about $p$-Laplacian equations; see
\cite{GK16,CFZ23,GL17,LGL15,LLF24,LFL14,SLW22,WJ24,WZ23,WW20,YCC23,YZ24} and their references.
Lattice dynamical systems are developed by spatial discretizations of
general partial differential equations,
where a countable system of ordinary differential equations is constructed
by replacing the spatial derivatives by differences (\cite{D77}).
Lattice dynamical systems have been applied in many important areas and attracted
much attention these years;
see \cite{B21,CKV14,CFZ23,GK16,GL17,HLW21,HSZ11,LLW22,LLW23,N20,SL22,WJ24,WWHK,XZL20}
and their references.
In particular, Wang and Jin recently considered the stochastic $p$-Laplacian lattice equation
with multiplicative noise in \cite{WJ24}, in which they obtained the existence of
pullback random $(\ell^2,\ell^q)$-attractors for each $q\in[1,+\8)$ and
constructed the invariant sample measures in $\ell^2$ for the
nonautonomous random dynamical system (NRDS for short) generated by \eqref{1.1}.

In the present article, based on the results of \cite{WJ24},
we mainly consider the convergence properties for the stochastic $p$-Laplacian
lattice equation \eqref{1.1} in the following four aspects:
\benu\item pullback random $(\ell^2,\ell^q)$-attractors for all $q\in[1,\8)$,
\item time-dependent invariant sample measures,
\item stochastic Liouville type equation, and
\item sample statistical solutions.\eenu

By convergence of pullback random $(\ell^2,\ell^q)$-attractors
$\cA_{\al}=\{\cA_{\al}(\tau,\w)\}_{\tau,\w}\in\R\X\W$
(see explanations in Subsection \ref{ss3.1}),
we mean the upper semi-continuity of $\cA_\al(\tau,\w)$ in $\ell^2\cap\ell^q$
as $\al\ra\al_0$ in $\R$, for each $\tau\in\R$ and $\bP$-a.s. $\w\in\W$,
where $(\W,\cF,\bP,\{\theta_t\}_{t\in\R})$ is an ergodic metric dynamical
system defined in Subsection \ref{ss2.1}.
The upper semi-continuity of attractors is a classical topic about attractors.
Upper semi-continuity of random single-space attractors of
stochastic partial differential equations
was considered in \cite{LCL14,LGL15,WLYJ21,XC22,ZZ23};
the case of random bi-spatial attractors was proved in \cite{LGL15}.
Wang studied the upper semi-continuity of single-space pullback random attractor
for the NRDS in \cite{W14}.
In \cite{LLF24} the authors obtained the upper semi-continuity of
tempered random attractors for stochastic delay $p$-Laplacian equation
on unbounded thin domains.
A result of upper semi-continuity of random bi-spatial pullback random attractors
of NRDS's was proved in \cite{CLY15} and it was applied to stochastic
parabolic equation on thin domains in \cite{LLW19} when the thin domain
collapses onto a lower dimensional domain.

But up to our knowledge until now, there is still no work about the upper
semi-continuity of pullback random $(\ell^2,\ell^q)$-attractors
for \eqref{1.1} for each $q\in[1,\8)$.
In our present work, we make some adjustment of Theorem 3.4 of \cite{CLY15}
to facilitate our estimates and show this convergence result in $\ell^2\cap\ell^q$.

For the convergence of invariant sample measures, we recall the theorem
of limiting behaviors for a sequence of invariant sample measures proved
by Chen and Yang in \cite{CY23}.
Originally, invariant sample measures are developed recently in \cite{CY23,WZ23,ZWC22}
for (nonautonomous) random dynamical systems from time-dependent invariant measures
(\cite{FMRT01,HLW21,LLW22,LLW23,WZC20,WZZ21,WWL22,YCC23}) for deterministic dynamical systems.
This sort of invariant measure is a ``weaker" version of the usual one for
random dynamical systems,
since the "sample" here represents the reliance of the invariant measures on each fixed $\w\in\W$,
while the usual one is in the sense of expectation in $\W$ (\cite{CFZ23,WLYJ21,WZL23,WW20,WCT23}).
Invariant sample measures has been obtained in \cite{WJ24}
for the stochastic $p$-Laplacian lattice equation \eqref{1.1} in $\ell^2$.
Then Theorem 4.1 of \cite{CY23} provides us an efficient method to obtain a weakly convergent
subsequence of a given sequence of families of invariant sample measures
constructed for \eqref{1.1} with $\al=\al_n$ and $\al_n\ra\al_0$ in $\R$.

Stochastic Liouville type equation we adopted in this article was developed in \cite{WZL24}
and originally established in \cite{CY23} by introducing It\^o's formula.
For Liouville type equation associated with deterministic differential systems, we
refer the reader to \cite{FMRT01,LLW22,LLW23}.
Liouville type equation was indeed taken as the basis for the equation
for time-dependent statistical solutions (\cite{FMRT01}).
Recently, it is extended to random Liouville type
equation of (nonautonomous) random differential systems (see \cite{CY23,WJ24,WZ23,ZWC22}).
In \cite{WZL24}, the authors provided a more general definition of the class
of test functions satisfying four pivotal conditions,
and hence developed stochastic Liouville type equation (\cite{CY23})
into more situations.
We will use the framework set up in \cite{WZL24} to check that the invariant sample
measures $\{\mu^{\al}_{\tau,\w}\}_{(\tau,\w)\in\R\X\W}$ in $\ell^2$ constructed in \cite{WJ24}
satisfy the following stochastic Liouville type equation,
for $s,t\in\R$ with $s\leqslant t$,
\begin{align}&\int_{H}\Psi(u)\mu^{\al}_{t,\theta_{t-\tau}\w}(\di u)
-\int_{H}\Psi(u)\mu^{\al}_{s,\theta_{s-\tau}\w}(\di u)\notag\\
=&\int_s^t\int_H\langle\tilde{F}(\vsig,u),\Psi'(u)\rangle
\mu^{\al}_{\vsig,\theta_{\vsig-\tau}\w}(\di u)\di\vsig
+\al\int_s^t\int_H(u,\Psi'(u))\mu^{\al}_{\vsig,\theta_{\vsig-\tau}\w}(\di u)
\di\tilde{W}(\vsig)\notag\\
&+\frac{\al^2}2\int_s^t\int_H\Psi''(u)(u,u)
\mu^{\al}_{\vsig,\theta_{\vsig-\tau}\w}(\di u)\di\vsig,\label{1.4}
\end{align}
for each test function $\Psi\in\cT$ (see definition in Subsection \ref{ss5.1}),
where $\tilde{F}$ is defined in \eqref{5.0A}
and $\tilde{W}(\vsig)=W(-\tau+\vsig)-W(-\tau)$.
In this part, the compactness of the $(\ell^2,\ell^q)$-attractor in $\ell^q$
plays an important role.

Note that the stochastic Liouville type equation depends on the parameter $\al$.
The weak convergence of invariant sample measures $\mu^{\al}_{\tau,\w}$
also implies convergence of each integral in the left side of the equation \eqref{1.4}.
This actually has asserted the convergence of the stochastic Liouville type equation
\eqref{1.4} in a coarse sense.
In fact, however, the termwise convergence of \eqref{1.4} as $\al$ varies is more interesting
and deserves more attention.
We will thereby consider and prove the termwise convergence of
the stochastic Liouville type equations.

Let $\{\mu^{\al}_{\tau,\w}\}_{(\tau,\w)\in\R\X\W}$ be the invariant sample measures constructed
in \cite{WJ24} for the NRDS generated by \eqref{1.1} for each $\al\in\R$.
In order to prove the termwise convergence of the stochastic Liouville type equation
with respect to $\al$,
one needs to further guarantee the convergence in $L^1$ and of the stochastic integrals
of the $\mu^{\al}_{\vsig,\theta_{\vsig-\tau}\w}$-integrals.
The convergence in $L^1$ follows from Lebesgue's Dominated Convergence Theorem,
and indeed from continuity of the solution $(t,\phi)\mapsto u_{\al}(t,\tau,\w,\phi)$
and compactness of the pullback random attractors.
Comparatively, the convergence of the stochastic integrals is not a natural consequence.
However, the convergence of other terms in the (stochastic Liouville type) equation
can well assure the convergence of the stochastic integral terms.

Statistical solutions were defined to describe the time evolution of
the probability distribution functions associated with the fluid flow (\cite{FMRT01}).
A significant feature of a statistical solution is its basis,
the Liouville type equation discussed above.
In \cite{CY23}, the authors defined a stochastic Liouville type equation for
a stochastic differential equation with respect to the invariant sample measures.
Afterwards, Wang, Zhu and Li not only developed the stochastic Liouville type equation,
but also gave a reasonable definition of \emph{sample statistical solution}
corresponding to the stochastic Liouville type equation in \cite{WZL24}.
In this article, we follow the definition of sample statistical solution
and prove that each family of invariant sample measures is
a sample statistical solution.
And then the convergence result of sample statistical solutions can be derived
from the convergence properties of invariant sample measures
stochastic Liouville type equations with a supplementary argument for the convergence
in $\cC([s,t])$ with $s\leqslant t$ by Arzel\`a-Ascoli theorem.

In the remainder of this article, Section \ref{s2} is about
the basic settings for this article and the NDRS induced by
\eqref{1.1}-\eqref{1.2}.
Section \ref{s3} is devoted to the convergence of $(\ell^2,\ell^q)$-pullback
random attractors for each $q\in[1,+\8)$.
In Section \ref{s4}, the convergence theorem of the invariant sample measures
is proved.
In Section \ref{s5}, we prove that the invariant sample measures satisfy
the stochastic Liouville type equation and obtain the termwise convergence of
the stochastic Liouville type equations.
In Section \ref{s6}, we prove that each family of invariant sample measures
is a statistical solution, which fulfills a convergence result as well.

\section{The induced nonautonomous random dynamical system}\label{s2}
	
	We introduce the basic settings of spaces and operators
	and induce the nonautonomous random dynamical system for \eqref{1.1}-\eqref{1.2}.
\subsection{Basic notations and properties}\label{ss2.1}
	For metric spaces $X$ and $Y$,
	we conventionally denote $\cC(X,Y)$ ($\cC_{\rm b}(X,Y)$) the collection of continuous
	(and bounded) functionals from $X$ to $Y$.
	When $Y=\R$, we simply use $\cC(X)$ ($\cC_{\rm b}(X)$) to represent $\cC(X,\R)$
	($\cC_{\rm b}(X,\R)$).
	
	Let $\ell^p$, $p\in[1,\8)$, be the Banach space of all real valued $p$-power summable
	bi-infinite sequences with the norm $\|\cdot\|_p$ such that
	$$\|u\|_p^p=\sum_{i\in\Z}|u_i|^p,\Hs\mb{for }u=(u_i)_{i\in\Z}\in\ell^p.$$
	For the case when $p=\8$, we let $\ell^\8$ be the Banach space of all bounded bi-infinite
	sequences with the norm $\|\cdot\|_\8$ satisfying
	$$\|u\|_\8=\sup_{i\in\Z}|u_i|,\Hs\mb{for }u=(u_i)_{i\in\Z}\in\ell^\8.$$
	For the case when $p=2$, $\ell^2$ is a Hilbert space with the inner product
	$$(u,v):=\sum_{i\in\Z}u_iv_i\hs\mb{for all }
	u=(u_i)_{i\in\Z},\,v=(v_i)_{i\in\Z}\in\ell^2$$
	and the norm simply denoted by $\|\cdot\|$.
Generally, we conventionally denote
$$(u,v):=\sum_{i\in\Z}u_iv_i\hs\mb{for all }u=(u_i)_{i\in\Z},\,v=(v_i)_{i\in\Z}$$
as long as the series is convergent.

	We have the following embedding
	relationship between $\ell^p$ and $\ell^q$, (see \cite[Lemma 2.1]{GL17} or \cite{WJ24})
	\be\label{2.1}1\leqslant p\leqslant q\leqslant\8\hs\Ra\hs \ell^p\subset \ell^q\mb{ and }
	\|u\|_q\leqslant\|u\|_p\mb{ for }u\in\ell^p.\ee
	For notational convenience in the sequel,
	given each $u\in\ell^p$, $p\geqslant1$ and $r>0$,
	we use $|u|$, $u^r$ and ${\rm sgn}u$ to denote the sequences as
	$$(|u|)_i=|u_i|,\Hs(u^r)_i=u_i^r\Hs\mb{and}\Hs({\rm sgn}u)_i={\rm sgn}u_i,$$
	respectively, for each $i\in\Z$,
	where ${\rm sgn}:\R\ra\{\pm1,0\}$ is the sign function, i.e.,
	${\rm sgn}x=x/|x|$ for $x\neq0$ and ${\rm sgn}0=0$.
	According to this, it is easy to see that $|u|^r=(|u_i|^r)_{i\in\Z}$
	for $u\in\ell^p$.
	We also set $0^0$ to be $1$ in the sequel.
	Let $u\in\ell^p$, $v\in\ell^q$.
	By \eqref{2.1} we have the following inequality:
	\be\label{2.2}(u,|v|^r)=\sum_{i\in\Z}u_i|v_i|^r\leqslant\sum_{i\in\Z}|u_i|\|v\|_{\8}^r
\leqslant\|u\|_1\|v\|_{q}^r\hs\mb{for }p=1,\,q\geqslant1,\,r\geqslant0.\ee
	
	Define a certain multiplication $\otimes$ for two arbitrary bi-infinite sequences such that
	$$(u\otimes v)_i:=u_iv_i\hs\mb{for all }
	u=(u_i)_{i\in\Z},\,v=(v_i)_{i\in\Z}.$$
	Define also the self-mappings $B$ and $B^*$ on the spaces of
bi-infinite sequences respectively as
	$$(Bu)_i:=u_{i+1}-u_i,\Hs(B^*u)_i:=u_{i-1}-u_i\Hs
\mb{with}\Hs Au:=B^*(|Bu|^{p-2}\otimes(Bu))$$
	$$\mb{i.e.,}\Hs
	(Au)_i=|u_i-u_{i-1}|^{p-2}(u_i-u_{i-1})-|u_{i+1}-u_i|^{p-2}(u_{i+1}-u_i).
	$$
	It is obvious that $(Bu,v)=(u,B^*v)$ for all
	$u=(u_i)_{i\in\Z},\,v=(v_i)_{i\in\Z}\in\ell^2$,
	both $B$ and $B^*$ are linear and bounded and
\be\label{2.3}\|Au\|^2\leqslant4\sum_{i\in\Z}|(Bu)_i|^{2p-2}\leqslant
2^{2p-1}\sum_{i\in\Z}\(|u_{i+1}|^{2p-2}+|v_i|^{2p-2}\)
=4^p\|u\|_{2p-2}^{2p-2}\leqslant4^p\|u\|^{2p-2},\ee
$$(Au,u)=\|Au\|\|u\|=2^p\|u\|^{p-1}\|u\|\leqslant2^p\|u\|^p.$$
	Then \eqref{1.1} can be rewritten as
	\be\label{2.4}\di u(t)=[-\nu(t)Au+f(t,u)-\lam u]\di t+\alpha u\di W(t),\; t>\tau,
	\hs u(\tau)=u_\tau\in\ell^2,\ee
	where $f:\R\X\ell^2\ra\ell^2$ satisfies that $(f(t,u))_i:=f_i(t,u_i)$,
	for each $i\in\Z$.

The Wiener process $W$ is determined as follows.
	Consider the space $\cC_0(\R)=\{\w\in\cC(\R):\w(0)=0\}$,
	with the compact open topology, the Borel $\sig$-algebra $\cF$ and
	the corresponding Wiener measure $\bP$.
	For sake of seeking for a solution of \eqref{2.4},
	we need to introduce the Ornstein-Uhlenbeck (OU for short) equation (see \cite{WLYJ21}).
	Let $\theta_t\w(\cdot)=\w(\cdot+t)-\w(t)$, $t\in\R$.
	Then $(\cC_0(\R),\cF,\bP,\{\theta_t\}_{t\in\R})$ is
	an ergodic metric dynamical system.
	Based on the system $(\cC_0(\R),\cF,\bP,\{\theta_t\}_{t\in\R})$, we set
	the stochastic stationary solution
	$$z(\theta_t\w)=-\int_{-\8}^0\me^{s}\theta_t\w(s)\di s,$$
	as a pathwise solution of the OU equation
	$\di z+z\di t=\di\w(t)$.
	Then there is a $\theta_t$-invariant set $\W\subset\cC_0(\R)$ such that
	$\bP(\W)=1$ and for every $\w\in\W$,
	$z(\theta_t\w)$ is continuous in $t$.
	The random variable $z(\theta_t\w)$ satisfies
	\be\label{2.5}
	\lim_{t\ra \pm\infty}\frac{z(\theta_t\w)}{t}
	=\lim_{t\ra \pm\infty}\frac{\int_{0}^{t}z(\theta_s\w)\di s}{t}=0.
	\ee
	Hence we only consider the set $\W$,
	and define the Wiener process $W$ on $(\W,\cF,\bP)$ as
	$$W(\w)(t):=\w(t)\Hs\mb{on the probability space }(\W,\cF,\bP).$$

\subsection{The nonautonomous random dynamical system}

We now present the already-known results in \cite{WJ24} here for the following discussion.
Let $u(t,\tau,\w,u_\tau)$ be a solution of \eqref{2.4} and set
\be\label{2.6}v(t)=\me^{-\al z(\theta_t\w)}u(t),\Hs
v(\tau)=v_\tau:=\me^{-\al z(\theta_\tau\w)}u_\tau.\ee
We deduce from \eqref{2.4} and \eqref{2.6} that
$v(t)$ satisfies the following random system
\be\label{2.7}
\frac{\di v(t)}{\di t}=F(t,\w,v),\hs t>\tau,\Hs v(\tau)=v_\tau,\ee
where
$$F(t,\w,v):=-\me^{\al(p-2)z(\theta_t\w)}\nu(t)Av+\(\al z(\theta_t\w)-\lam\)v
+\me^{-\al z(\theta_t\w)}f(t,\me^{\al z(\theta_t\w)}v).$$
	
We present the following assumption on the nonautonomous nonlinear forcing $f$ below.
	
\noindent\textbf{(F1)} The mapping $(t,u_i)\mapsto f_i(t,u_i)$ is continuous
	in both variables and differentiable in $u_i$,
	and there are $\beta>0$, $q\geqslant1$ and $\psi_k(t)=(\psi_{ki}(t))_{i\in\Z}$, $k=1,\,2$,
such that
	\begin{align}\label{2.8}f_i(t,u_i)u_i&\leqslant-\beta|u_i|^q+\psi_{1i}(t),
		\mb{ with }\psi_1\in L_{\rm loc}^1(\R;\ell^1)\\
\mb{and}\hs\frac{\pa f_i}{\pa u_i}(t,u_i)&\leqslant\psi_{2i}(t),
		\mb{ with }\psi_2\in L_{\rm loc}^{\8}(\R;\ell^{\8}).\notag\end{align}
	
Under the assumption \textbf{(F1)}, the equation \eqref{2.7} possesses a unique global solution
$v(\cdot,\tau,\w,v_\tau)\in\cC^1([\tau,\8);\ell^2)\cap L^q_{\rm loc}(0,\8;\ell^q)$
for each $(\tau,v_\tau)\in\R\X\ell^2$ and $\bP$-a.s. $\w\in\W$,
such that $v(\tau,\tau,\w,v_\tau)=v_\tau$,
where $\cC^1$ means differentiability and continuity of the derivative.
Then set
	\be\label{2.9}u(t,\tau,\w,u_\tau)=\me^{\al z(\theta_t\w)}
v(t,\tau,\w,\me^{-\al z(\theta_\tau\w)}u_\tau).\ee
	We know that $u(t,\tau,\w,u_\tau)$ is exactly the unique global solution of \eqref{2.4}
	for each $(\tau,u_\tau)\in\R\X\ell^2$ and $\bP$-a.s. $\w\in\W$.
	
	We define a mapping $\vp:\R^+\X\R\X\W\X\ell^2\ra\ell^2$,
	such that for every $(t,\tau,\w,u_\tau)\in\R^+\X\R\X\W\X\ell^2$,
	\be\label{2.10}\vp(t,\tau,\w,u_\tau)=u(t+\tau,\tau,\theta_{-\tau}\w,u_\tau)
	=\me^{\al z(\theta_t\w)}v(t+\tau,\tau,\theta_{-\tau}\w,v_\tau),\ee
	with $v_\tau=\me^{-\al z(\theta_\tau\w)}u_\tau$.
	It has been proved in \cite{WJ24} that $\vp$ defined in \eqref{2.10} is actually
	a continuous NRDS over $(\W,\cF,\bP,\{\theta_t\}_{t\in\R})$,
	i.e., $\vp$ satisfies the following conditions:
	\benu\item[(1)] $\vp(\cdot,\tau,\cdot,\cdot):\R^+\X\W\X\ell^2\ra\ell^2$ is
	$(\cB(\R^+)\X\cF\X\cB(\ell^2),\cB(\ell^2))$-measurable for every $\tau\in\R$;
	\item[(2)] $\vp(0,\tau,\w,\cdot)$ is the identity on $\ell^2$ for every $(\tau,\w)\in\R\X\W$;
	\item[(3)] $\vp(t+s,\tau,\w,\cdot)=\vp(t,\tau+s,\theta_s\w,\vp(s,\tau,\w,\cdot))$
	for every $t,s\in\R^+$ and $(\tau,\w)\in\R\X\W$;
	\item[(4)]$\vp(\cdot,\tau,\w,\cdot):\R^+\X\ell^2\ra\ell^2$ is
continuous for every $(\tau,\w)\in\R\X\W$
	for $\bP$-a.s. $\w\in\W$.
	\eenu

\section{Convergence of pullback random attractors}\label{s3}
\subsection{Preliminaries on pullback random attractors}\label{ss3.1}
	
Let $X$ be a Polish space and $2^X$ be the collection of all subsets of $X$.
An $X$-valued time-parametrized random variable $R(\tau,\w)$ is called
a \emph{nonautonomous random variable} in $X$.
A \emph{nonautonomous random set} $D$ in $X$ is a family of nonempty
subsets $D(\tau,\w)\in 2^X$ with two parameters
$(\tau,\w)\in\R\X\W$, which is $(\W,\cF,\bP)$-measurable ($\cF$-measurable for short)
with respect to $\w\in\W$ (see \cite{C02,W12}),
i.e., the mapping $\w\mapsto \di_X(x,D(\tau,\w))$ is $(\cF,\cB(\R))$-measurable for each fixed
$(x,\tau)\in X\X\R$, where $\di_X(\cdot,\cdot)$ is the metric induced by norm of $X$.

In order to present the concept of bi-spatial pullback attractors for an NRDS,
we consider two separable Banach spaces $(X,\|\cdot\|_X)$ and $(Y,\|\cdot\|_Y)$.
An NRDS $\vp$ on $X$ is assumed to take its values in the terminate space
$Y$, i.e.,
$$\vp(t,\tau,\w,X)\subset Y,\;\mb{for all }t,\,\tau\in\R
\mb{ with }t>0\mb{ and }\bP\mb{-a.s. }\w\in\W.$$
For convergence uniformity, we assume the space pair $(X,Y)$ is \emph{limit-identical}
in the following sense,
$$
x_n\in X\cap Y,\hs
\|x_n-x_0\|_X\ra0 \mb{ and }\|x_n-y_0\|_Y\ra0
\mb{ imply }x_0=y_0\in X\cap Y.
$$
The norm of $X\cap Y$ is often defined as $\|x\|_{X\cap Y}=\|x\|_X+\|x\|_Y$ for all
$x\in X\cap Y$, which makes $(X\cap Y,\|\cdot\|_{X\cap Y})$ a normed linear space.
It was proved in \cite{LGL15} that, $(X\cap Y,\|\cdot\|_{X\cap Y})$ is a Banach space
if and only if $(X,Y)$ is a limit-identical pair.
Moreover, when $(X,Y)$ is limit-identical and $A\subset X\cap Y$, $A$ is compact in
$X\cap Y$ if and only if $A$ is compact in $X$ and $Y$, respectively.

Let $\cA=\{\cA(\tau,\w)\}_{(\tau,\w)\in\R\X\W}$ be a nonautonomous random set in $X\cap Y$
and $D=\{D(\tau,\w)\}_{(\tau,\w)\in\R\X\W}$ be a nonautonomous random set in $X$.
We say that $A$ \emph{pullback attracts} in (the topology of) $X$ if
\be\label{3.1}\lim_{t\ra +\8}{\rm dist}_{X}(\vp(t,\tau-t,
\theta_{-t}\w,D(\tau-t,\theta_{-t}\w)),\cA(\tau,\w))=0,\ee
where ${\rm dist}_X(\cdot,\cdot)$ is the \emph{Hausdorff semi-distance}
under the norm of an arbitrary normed space $X$, i.e., for two nonempty sets $A,B\subset X$,
$${\rm dist}_X(A,B):=\sup_{a\in A}\inf_{b\in B}\|a-b\|_X.$$
Here the pullback attraction can also be taken in (the topology of) $X\cap Y$ instead
with $X$ in \eqref{3.1} replaced by $X\cap Y$.

We say $\cD$ is a \emph{universe} in $X$, if $\cD$ is
an inclusion-closed collection of nonautonomous
random sets, i.e., when $D\in\cD$ and a nonautonomous random set $D'$ satisfies
$D'(\tau,\w)\subset D(\tau,\w)$ for each $(\tau,\w)\in\R\X\W$, then $D'\in\cD$.
	
In this article, we adopt the much weaker definitions of bi-spatial pullback (random) attractors
presented in \cite{CLL18,WJ24} as follows, where we do not require the attractor
to be $\cF$-measurable
in the initial phase space.
The existence theorem of pullback random bi-spatial attractors can be found in \cite{WJ24}.

\bd\label{de3.1}
Let $\vp$ be an NRDS on $X$ taking its values in $Y$ and
$\cD$ be a universe in $X$.
A set-mapping $\cA:\R\X\W\ra 2^{X\cap Y}$ is said to be a
\textbf{pullback random $(X,Y)$-attractor} for $\vp$ with respect to $\cD$ if
\benu\item[(1)] $\cA(\tau,\w)$ is compact in $Y$ for all $\tau\in\R$ and $\bP$-a.s. $\w\in\W$;
	\item[(2)] $\cA$ is invariant under the system $\vp$, i.e.,
	$$\vp(t,\tau,\w,\cA(\tau,\w))=\cA(\tau+t,\theta_t\w),
	\mb{ for all }t\geqslant0\mb{ and }(\tau,\w)\in\R\X\W;$$
	\item[(3)] $\cA$ is pullback $\cD$-attracting in $X\cap Y$, that is, $\cA$ pullback attracts
every $D\in\cD$ in the topology of $X\cap Y$,
		\item[(4)] there exists a closed set $K\in\cD$ with $\cA\subset K$.
	\item[(5)] $\cA$ is a nonautonomous random set in $Y$.
	\eenu
\ed

Let $\vp$ be an NRDS on $X$ taking its values in $Y$ and $\cD$ be a universe.
A \emph{pullback $\cD$-absorbing set} in $Y$ for the NRDS $\vp$
is a nonautonomous random closed set $K$ in $Y$, such that
for every $D\in\cD$ and all $\tau\in\R$, $\bP$-a.s. $\w\in\W$,
there exists $T=T(\tau,\w,D)>0$ such that for all $t>T$,
$$\vp(t,\tau-t,\theta_{-t}\w,D(\tau-t,\theta_{-t}\w))\subset K(\tau,\w).$$

For the upper semi-continuity of pullback random attractors,
we consider a family $\{\vp_\al\}_{\al\in\fA}$ of NRDS's with $(\fA,\rho)$ a metric space.
We can find the following theorem in \cite{CLY15} with a slight modification
(see more in \cite{LCL14,LGL15,W14}) for the upper semi-continuity.
	
\bt\label{th3.2}
Let $(X,Y)$ be a limit-identical pair of separable Banach spaces
and $\vp_\al$ be a continuous NRDS on $X$ over $(\W,\cF,\bP,\{\theta_t\}_{t\in\R})$
taking values in $Y$ for each $\al\in\fA$ with $\fA$ a metric space.
For each $\al\in\fA$, let $\cD_{\al}$ be a universe for $\vp_\al$ in $X$.
Suppose there exists $\al_0\in\fA$ such that
\benu\item[(1)] for every $t>0$, $\tau\in\R$, $\w\in\W$,
	$\al_n\in\fA$ with $\al_n\ra\al_0$ and $x_n,x\in X$ with $x_n\ra x$,
	\be\label{3.2}
	\lim_{n\ra\8}\vp_{\al_n}(t,\tau,\w,x_n)=\vp_{\al_0}(t,\tau,\w,x)\hs\mb{ in }X;
	\ee
\item[(2)] for each $\al\in\fA$, $\vp_\al$ has
a pullback $\cD_\al$-absorbing set $K_\al\in\cD_\al$ in $X$ such that
    for each $(\tau,\w)\in\R\X\W$ and convergent sequence $\al_n\ra\al_0$ such that
	\be\label{3.3}\Cap_{n\geqslant 1}\ol{\Cup_{N\geqslant n}
K_{\al_N}(\tau,\w)}^X\subset K_{\al_0}(\tau,\w)\ee
	and $K_{\al_0}=\{K_{\al_0}(\tau,\w)\}_{(\tau,\w)\in\R\X\W}\in\cD_{\al_0}$;
\item[(3)] for each $\al\in\fA$, $\vp_\al$ has
a pullback random $(X,Y)$-attractor $\cA_\al\in\cD_\al$
with respect to $\cD_\al$ and there exists a neighborhood $\fU$ of $\al_0$ in $\fA$ such that
for every $(\tau,\w)\in\R\X\W$,
\be\label{3.4}\Cup_{\al\in\fU}\cA_{\al}(\tau,\w)\hs\mb{is precompact in }X\cap Y.\ee
\eenu
Then for every $(\tau,\w)\in\R\X\W$,
\be\label{3.5}\lim_{\al\ra\al_0}{\rm dist}_{X\cap Y}
(\cA_\al(\tau,\w),\cA_{\al_0}(\tau,\w))=0.\ee
\et
\bo The proof is similar to that of Proposition 3.2 and Theorem 3.4 of \cite{CLY15}.
For the reader's convenience, we give a brief verification.
	
Suppose the contrary. We thus have $(\tau,\w)\in\R\X\W$, a positive number $a$
and a sequence $\al_n\ra\al_0$ in $\fU$ such that for all $n\in\N^+$,
${\rm dist}_{X\cap Y}(\cA_{\al_n}(\tau,\w),\cA_{\al_0}(\tau,\w))\geqslant2a$,
and hence the existence of a sequence $\{x_n\}_{n\in\N_+}$
with $x_n\in\cA_{\al_n}(\tau,\w)$ such that
\be\label{3.6}{\rm dist}_{X\cap Y}(x_n,\cA_{\al_0}(\tau,\w))
\geqslant a\hs\mb{for all }n\in\N_+.\ee
By \eqref{3.4}, we can as well assume (pick a subsequence if necessary) that
\be\label{3.7}\lim_{n\ra\8}x_n=x_0\hs\mb{in }X\cap Y.\ee
Then we claim $x_0\in\cA_{\al_0}(\tau,\w)$,
which contradicts \eqref{3.6} and completes the proof.
	
It suffices to show the claim.
Take a positive time sequence $\{t_m\}_{m\in\N_+}$ with $t_m\ra+\8$.
Fix $m=1$.
For each $n\in\N^+$, by the invariance of $\cA_{\al_n}$, there exists
$x_{1,n}\in\cA_{\al_n}(\tau-{t_1},\theta_{-t_1}\w)$ such that,
\be\label{3.8}x_n=\vp_{\al_n}(t_1,\tau-t_1,\theta_{-t_1}\w,x_{1,n}).\ee
Note that $x_{1,n}\in\cA_{\al_n}(\tau-t_1,\theta_{-t_1}\w)$.
By \eqref{3.4} again, we can similarly assume that $x_{1,n}\ra y_1$ in $X\cap Y$ as $n\ra\8$,
for some $y_1\in X\cap Y$.
Using \eqref{3.7}, \eqref{3.8} and \eqref{3.2}, we have
$$x_0=\lim_{n\ra\8}\vp_{\al_n}(t_1,\tau-t_1,\theta_{-t_1}\w,x_{1,n})
=\vp_{\al_0}(t_1,\tau-t_1,\theta_{-t_1}\w,y_1)\hs\mb{in }X.$$
Recalling that $\cA_{\al}\subset K_{\al}$ and $x_{1,n}\in\cA_{\al_n}(\tau-t_1,\theta_{-t_1}\w)$
for each $n\in\N^+$, we know by \eqref{3.3} that
$$y_1=\lim_{n\ra\8}x_{1,n}\in\Cap_{n\geqslant 1}\ol{\Cup_{N\geqslant n} K_{\al_N}(\tau,\w)}^X
\subset K_{\al_0}(\tau-t_1,\theta_{-t_1}\w).$$
By induction, for each $m\in\N^+$, we can find
a $y_m\in K_{\al_0}(\tau-t_m,\theta_{-t_m}\w)$ such that
$$x_0=\vp_{\al_0}(t_m,\tau-t_m,\theta_{-t_m}\w,y_m).$$
Using the attraction property of $\cA_{\al_0}$ in $\cD_{\lam_0}$, we obtain
\begin{align*}&{\rm dist}_{X}(x_0,\cA_{\al_0}(\tau,\w))
={\rm dist}_{X}(\vp_{\al_0}(t_m,\tau-t_m,\theta_{-t_m}\w,y_m),\cA_{\al_0}(\tau,\w))\\
\leqslant&{\rm dist}_{X}(\vp_{\al_0}(t_m,\tau-t_m,
\theta_{-t_m}\w,K_{\al_0}(\tau-t_m,\theta_{-t_m}\w)),\cA_{\al_0}(\tau,\w))\ra0,
\end{align*}
as $m\ra\8$, which implies $x_0\in\cA_{\al_0}(\tau,\w)$ by its compactness in $Y$
and proves the claim.
\eo
\br If we replace $X\cap Y$ by $Y$ in \eqref{3.4} and \eqref{3.5},
Theorem \ref{th3.2} would be a simple deduction of \cite[Theorem 3.4]{CLY15}.
In Theorem \ref{th3.2}, we obtain a stronger conclusion by a stronger condition.
\er
The condition (3) of Theorem \ref{th3.2} can be also inferred
by the local uniform pullback asymptotic compactness of $\vp_\al$
(see \cite{LCL14,LGL15} for the case without universe).
\bd\label{de3.3} Let $\{\vp_\al\}_{\al\in\fA}$ be a family of NRDS's and
$\cD_{\al}$ be a universe of $\vp_\al$ in $X$ for each $\al\in\fA$.
We say that $\{\vp_\al\}_{\al\in\fA}$ is
\textbf{locally uniformly pullback $\cD_\al$-asymptotically compact}
in $X$ (resp. in $Y$), if for each $\al_0\in\fA$,
there exists a neighborhood $\fU$ of $\al_0$ in $\fA$
such that for each $(\tau,\w)\in\R\X\W$,
whenever there are sequences $\al_n\in\fU$, $0<t_n\ra+\8$ and
$x_n\in D_{\al_n}(\tau-t_n,\theta_{-t_n}\w)$
with $D_{\al_n}\in\cD_{\al_n}$, the sequence
$\{\vp_{\al_n}(t_n,\tau-t_n,\theta_{-t_n}\w,x_n)\}_{n\in\N^+}$
has a convergent subsequence in $X$ (resp. in $Y$).
If $\al_n\equiv\al$ for all $n\in\N^+$,
we say $\vp_\al$ is \textbf{pullback $\cD_\al$-asymptotically compact} in $X$ (resp. in $Y$).\ed

\bt\label{th3.4A} Let $\{\vp_\al\}_{\al\in\fA}$ be a family of NRDS's with
$\cD_{\al}$ a universe of $\vp_\al$ in $X$ for each $\al\in\fA$ and $\al_0\in\fA$.
Suppose that
\benu\item[(3$'$)]$\{\vp_\al\}_{\al\in\fA}$ is locally uniformly
pullback $\cD_\al$-asymptotically compact in $X$ and $Y$,
respectively and each $\vp_\al$ has a pullback random $(X,Y)$-attractor
$\cA_\al\in\cD_\al$ with respect to $\cD_\al$.
\eenu
Then there is a neighborhood $\fU$ of $\al_0$ in $\fA$ such that for every $(\tau,\w)\in\R\X\W$,
the set $\Cup_{\al\in\fU}\cA_{\al}(\tau,\w)$ is precompact in $X\cap Y$.
As a result, the conditions (1), (2) and (3$'$) can also ensure
the upper semi-continuity \eqref{3.5}.
\et
\bo By (3$'$) and Definition \ref{de3.3}, we let $\fU$ be a neighborhood of $\al_0$ in $\fA$,
such that
for arbitrarily fixed $(\tau,\w)\in\R\X\W$,
whenever there are sequences $\al_n\in\fU$, $0<t_n\ra+\8$ and
$x_n\in D_{\al_n}(\tau-t_n,\theta_{-t_n}\w)$ with $D_{\al_n}\in\cD_{\al_n}$,
the sequence $\{\vp_{\al_n}(t_n,\tau-t_n,\theta_{-t_n}\w,x_n)\}_{n\in\N^+}$
has a convergent subsequence in $X$ and $Y$, respectively.

In order to show the precompactness of
$$\sA(\tau,\w):=\Cup_{\al\in\fU}\cA_{\al}(\tau,\w)$$
in $X$, we pick a sequence $\{y_n\}$ in $\sA$.
Then by the invariance of attractors,
we can find $\al_n\in\fU$ and $x'_n\in\ell^2$ such that $y_n\in\cA_{\al_n}(\tau,\w)$,
$x'_n\in\cA(\tau-n,\theta_{-n}\w)$ and
$$\vp_{\al_n}(n,\tau-n,\theta_{-n}\w,x'_n)=y_n.$$
Noting that $\cA_{\al_n}\in\cD_{\al_n}$,
we actually have already ensured the precompactness of $\{y_n\}_{n\in\N^+}$ and
hence of $\sA(\tau,\w)$ in $X$
by the local uniform pullback $\cD_\al$-asymptotic compactness of $\{\vp_\al\}$ in $X$.

The conclusion in $Y$ can be obtained by a similar argument in $Y$.
Since $(X,Y)$ is a limit-identical pair, $\sA(\tau,\w)$ is obviously precompact in $X\cap Y$.
The proof is finished.
\eo
\subsection{Convergence of pullback random bi-spatial attractors}

In order to emphasize the coefficient $\al$ in \eqref{2.4},
we always denote the solutions of \eqref{2.4} and \eqref{2.7} by $u_\al$ and $v_\al$, respectively
and by $\vp_\al$ the generated process by $u_\al$ in the sequel.

Let $\lam_0,\,\lam_1\in\R$ be arbitrary such that
$$0<\lam_1<\lam_0<2\lam.$$
We recall some other assumptions in \cite{WJ24} on $f$ with a slight adjustment for the following argument:

\noindent\textbf{(F2)} The function $\psi_1$ in \eqref{2.8} satisfies
$$\int_{-\8}^0\me^{\lam_1s}\|\psi_1(s)\|_1\di s<+\8.$$

\noindent\textbf{(F3)} Let $q\in[1,+\8)$ be given in \textbf{(F1)}.
Assume that there are $\psi_k(t)=(\psi_{ki}(t))_{i\in\Z}$, $k=3,4$, such that
\be\label{3.1A}|f_i(t,u_i)|\leqslant\psi_{3i}(t)|u_i|^{q-1}+\psi_{4i}(t),\ee
$$\mb{with}\hs\psi_3\in L^q_{\rm loc}(\R,\ell^q)\cap L_{\rm loc}^1(\R,\ell^1),
\hs\psi_4\in L_{\rm loc}^1(\R,\ell^1),$$
and if moreover $q\in[1,2)$, there are $\kappa_0$, $\Lam$, $t_0>0$
with $q\lam_0>2\kappa_0$ such that when $t>t_0$,
$$\|\psi_3(-t)\|_1\leqslant\Lam\me^{\kappa_0t}\hs\mb{and}
\hs\|\psi_4(-t)\|_1\leqslant\Lam\me^{\kappa_0t}.$$

According to the settings in \cite{WJ24}, for each $\vp_\al$ ($\al\in\R$),
we consider the universe $\cD_\al$ to be all the families
$D_\al=\{D_\al(\tau,\w)\}_{(\tau,\w)\in\R\X\W}$
be a nonautonomous random set of $\ell^2$ such that for every $\tau\in\T$ and $\w\in\W$,
$$\lim_{s\ra-\8}\me^{\lam_0s+2\al\int_s^0
z(\theta_r\w)\di r-2\al z(\theta_t\w)}\|D_\al(\tau+s,\theta_s\w)\|^2=0,$$
where $\|D\|=\sup_{u\in D}\|u\|$ for each subset $D$ of $\ell^2$.
The following theorem is a main theorem --- Theorem 3.15 of \cite{WJ24}.

\bt\label{th3.5} Let the assumptions \textbf{(F1)}, \textbf{(F2)} and \textbf{(F3)} hold for $f$.
Then the NRDS $\vp_{\al}$ has a unique pullback
random $(\ell^2,\ell^q)$-attractor $\cA_\al\in\cD_\al$
with respect to $\cD_\al$,
which is also the pullback random $(\ell^2,\ell^2)$-attractor with respect to $\cD_\al$.
\et
In the following discussion, we are devoted to verify the upper
semi-continuity of $\cA_\al$ at each $\al\in\R$.
By Theorems 2.2 of \cite{WJ24}, we have the following consequence.
\bl For each $u_\tau\in\ell^2$, $\al\in\R$, $\tau\in\R$, $T>\tau$,
$t\in(\tau,T]$ and a.s. $\w\in\W$,
it holds that
\be\label{3.11}
\|v_{\al}(t)\|^2\leqslant G_1(\al,\|u_\tau\|)\Hs\mb{and}\hs
\int_\tau^t\|v_\al(s)\|_q^q\di s\leqslant G_2(\al,\|u_\tau\|)\ee
where for each $R\in\R^+$,
\begin{align*}
&G_1(\al,R):=G_1(\al,\tau,T,\w,R)\\
=&\me^{2\al\int_\tau^T|z(\theta_r\w)|\di r-2\al z(\theta_\tau\w)}R^2
+2\int_\tau^T\me^{2\al\int_s^T|z(\theta_r\w)|\di r-2\al z(\theta_s\w)}\|\psi_1(s)\|_1\di s,\\
&G_2(\al,R):=G_2(\al,\tau,T,\w,R)=E(\al)G_1(\al,R)\end{align*}
and
$$E(\al):=E(\al,\tau,T,\w)=\frac{1}{2\beta}\me^{2\lam(T-\tau)+2|\al|\int_\tau^T|z(\theta_r\w)|\di r
+|\al(q-2)|\max_{s\in[\tau,T]}|z(\theta_s\w)|}.$$
Furthermore, $G_i(\al,R)$ ($i=1,2$) is joint continuous in $(\al,R)\in\R\X\R^+$.
\el
\bo The proof is almost the same as that of Theorems 2.2 of \cite{WJ24}.
We only give a simple sketch of the proof for the reader's convenience.

Taking the inner product of \eqref{2.7} and $v$, and applying Gronwall's lemma
and the assumption (F1), we can easily obtain (recall the inequality (18) of \cite{WJ24})
\begin{align}
&\|v_{\al}(t)\|^2+2\beta\int_\tau^t\me^{2\lam(\vsig-t)+2\al\int_\vsig^tz(\theta_r\w)\di r
+\al(q-2)z(\theta_\vsig\w)}\|v_\al(\vsig)\|_q^q\di\vsig\notag\\
\leqslant&
\me^{2\lam(\tau-t)+2\al\int_\tau^tz(\theta_r\w)\di r}\|v_\al(\tau)\|^2
+2\int_\tau^t\me^{2\lam(s-t)+2\al\int_s^tz(\theta_r\w)\di r-2\al z(\theta_s\w)}\|\psi_1(s)\|_1\di s,
\label{3.12}\end{align}
from which the first inequality of \eqref{3.11} follows by \eqref{2.6} as well.

It is trivial to see that
$$\frac{1}{E(\al)}\int_\tau^T\|v_\al(s)\|^2\di s\leqslant
2\beta\int_\tau^t\me^{2\lam(\vsig-t)+2\al\int_\vsig^tz(\theta_r\w)\di r
+\al(q-2)z(\theta_\vsig\w)}\|v_\al(\vsig)\|_q^q\di\vsig\leqslant G_1(\al,\|u_\tau\|),$$
which implies the second inequality of \eqref{3.11}.

The joint continuity of $G_i(\al,R)$ ($i=1,2$) in $(\al,R)$ is
a simple deduction of the Dominated Convergence Theorem.
The proof is complete.
\eo	
In the following argument, for notational convenience, we set
$$\cZ=\(1\vee\sup_{n\in\N}|\al_n|\)\cdot\max_{t\in[\tau,T]}|z(\theta_t\w)|,\hs\hs
\cE_{r,n}:=\max_{t\in[\tau,T]}\left|\me^{\alpha_nrz(\theta_t\w)}
-\me^{\alpha_0rz(\theta_t\w)}\right|$$
$$\mb{and}\hs R_\de:=R_\de(u_\tau)=\sup\{\|u\|:\|u-u_\tau\|<\de,\,u\in\ell^2\}
\hs\mb{for each }\de>0,$$
where $a\vee b$ means the bigger one between the real numbers $a$ and $b$, $r\in\R$.
It is obvious that for each $r\in\R$ and $\w\in\W$, $\cE_{r,n}\ra0$ as $n\ra\8$.
By continuity of $G_1$ given in Lemma \ref{th3.2},
we know that for each $t\in[\tau,T]$, $\w\in\W$ and $u'_\tau$
with $\|u'_\tau-u_\tau\|<\de$,
$$\|v_{\al_n}(t,\tau,\w,u'_\tau)\|^2\leqslant\sup_{n\in\N}G_1(\al_n,R_\de):=\cG_\de<\8.$$
All $\cZ$, $\cE_{r,n}$ and $\cG_\de$ are independent of the time $t\in[\tau,T]$.
We fix an arbitrary $\de>0$ and consider the initial $u'_\tau$ with $\|u'_\tau-u_\tau\|<\de$.
Moreover, we always denote $C$ as an arbitrary positive constant only depending on $p$ and $q$
and $C$ may be different from line to line and even in the same line.

Next we check the condition (1) of Theorem \ref{th3.2} in the following lemma.
\bl\label{le3.7} Suppose the assumptions \textbf{(F1)}, \textbf{(F2)} and \textbf{(F3)} hold.
Let $\al_0\in\R$ with a sequence $\{\al_n\}\subset\R$ such that $\al_n\ra\al_0$ as $n\ra\8$,
$\tau,\,t,\,T\in\R$ with $\tau<t\leqslant T$, $u_\tau\in\ell^2$ and $\bP$-a.s. $\w\in\W$.
Then there exists $\de:=\de(\tau,T)>0$ such that when $\|\phi-u_\tau\|<\de$,
$$u_{\al_n}(t,\tau,\w,\phi)\ra u_{\al_0}(t,\tau,\w,\phi)$$
in $\ell^2$ uniformly in $\phi$ with $\|\phi-u_\tau\|<\de$ as $\al_n\ra\al_0$.\el
\bo It suffices to prove
\be\label{3.15}v_{\al_n}(t,\tau,\w,\me^{-\al_n z(\theta_\tau\w)}\phi)\ra
v_{\al_0}(t,\tau,\w,\me^{-\al_0 z(\theta_\tau\w)}\phi)\ee
in $\ell^2$ uniformly in $\phi$ with $\|\phi-u_\tau\|<\de$ as $\al_n\ra\al_0$
for some $\de>0$.
This is due to the fact
\begin{align*}&\|u_{\al_n}(t,\tau,\w,\phi)-u_{\al_0}(t,\tau,\w,\phi)\|\\
\leqslant&\me^{\max\{|\al_n|:n\in\N\}|z(\theta_t\w)|}\|v_{\al_n}(t)-v_{\al_0}(t)\|
+|\me^{\al_nz(\theta_t\w)}-\me^{\al_0 z(\theta_t\w)}|\|v_{\al_0}(t)\|,\end{align*}
which obviously vanishes as $\al_n\ra\al_0$.

In order to show \eqref{3.15}, we take the inner product in $\ell^2$ of the difference
between the corresponding equations \eqref{2.7} with $\al_n$ and $\al_0$ and $v_{\al_n}-v_{\al_0}$.
Then we have
\begin{align}
&\frac{1}{2}\frac{\di}{\di t}\|v_{\al_n}-v_{\al_0}\|^2=
\(F(t,\w,v_{\al_n})-F(t,\w,v_{\al_0}),v_{\al_n}-v_{\al_0}\)\notag\\
=&-\nu(t)\(\me^{\alpha_n(p-2)z(\theta_t\w)}Av_{\al_n}
-\me^{\alpha_0(p-2)z(\theta_t\w)}Av_{\al_0},v_{\alpha_n}-v_{\alpha_0}\)\notag\\
&+\((\al_nz(\theta_t\w)-\lam)v_{\al_n}
-(\al_0 z(\theta_t\w)-\lam)v_{\al_0},v_{\al_n}-v_{\al_0}\)\notag\\
&+\(\me^{-\alpha_nz(\theta_t\w)}f(t,\me^{\alpha_nz(\theta_t\w)}v_{\al_n})
-\me^{-\alpha_0 z(\theta_t\w)}f(t,\me^{\alpha_0 z(\theta_t\w)}v_{\al_0}),v_{\al_n}
-v_{\al_0}\)\notag\\
:=&\nu(t)J_1+J_2+J_3.\label{3.16}
\end{align}
First, observe by \eqref{2.3}, \eqref{2.2} and the mean value theorem that
\begin{align}|J_1|\leqslant&
\left|\((\me^{\alpha_n(p-2)z(\theta_t\w)}-\me^{\alpha_0(p-2)z(\theta_t\w)})
Av_{\al_n},v_{\alpha_n}-v_{\alpha_0}\)\right|\notag\\
&+\me^{\alpha_0(p-2)z(\theta_t\w)}\left|\((Av_{\al_n}-Av_{\alpha_0}),
v_{\alpha_n}-v_{\alpha_0}\)\right|\notag\\
:=&J_{11}+J_{12},\\
J_{11}\leqslant&\cE_{p-2,n}\|Av_{\al_n}\|(\|v_{\alpha_n}\|+\|v_{\alpha_0}\|)
\leqslant C\cG_\de^{\frac{p}2}\cE_{p-2,n},\end{align}
\begin{align}J_{12}\leqslant&\me^{(p-2)\cZ}\sum_{i\in\Z}\left||(Bv_{\al_n})_i|^{p-2}
(Bv_{\al_n})_i-|(Bv_{\al_0})_i|^{p-2}(Bv_{\al_0})_i\right|
|(Bv_{\alpha_n}-Bv_{\alpha_0})_i|\notag\\
\leqslant&(p-1)\me^{(p-2)\cZ}\sum_{i\in\Z}\(|(Bv_{\al_n})_i|
+|(Bv_{\al_0})_i|\)^{p-2}|(Bv_{\alpha_n}-Bv_{\alpha_0})_i|^2\notag\\
\leqslant&(p-1)\me^{(p-2)\cZ}\(\|Bv_{\al_n}\|+\|Bv_{\al_0}\|\)^{p-2}
\|B(v_{\alpha_n}-v_{\alpha_0})\|^2\notag\\
\leqslant&C\me^{(p-2)\cZ}\cG_{\de}^{\frac{p-2}{2}}\|v_{\alpha_n}-v_{\alpha_0}\|^2
\end{align}
and
\begin{align}
J_2=&(\al_n-\al_0)z(\theta_t\w)(v_{\al_n},v_{\al_n}-v_{\al_0})
+(\al_0z(\theta_t\w)-\lam)\|v_{\al_n}-v_{\al_0}\|^2\notag\\
\leqslant&C\cZ\cG_\de|\al_n-\al_0|+|\cZ-\lam|\|v_{\al_n}-v_{\al_0}\|^2.
\end{align}
We divide $J_3$ into three parts as follows:
\begin{align}|J_3|\leqslant
&\left|\me^{-\alpha_nz(\theta_t\w)}-\me^{-\alpha_0z(\theta_t\w)}\right|
\left|\(f(t,\me^{\alpha_nz(\theta_t\w)}v_{\al_n}),v_{\al_n}-v_{\al_0}\)\right|\notag\\
&+\me^{-2\alpha_0z(\theta_t\w)}\left|\(f(t,\me^{\alpha_nz(\theta_t\w)}v_{\al_n})
-f(t,\me^{\alpha_0z(\theta_t\w)}v_{\al_0}),
\me^{\al_0 z(\theta_t\w)}v_{\al_n}-\me^{\al_nz(\theta_t\w)}v_{\al_n}\)\right|\notag\\
&+\me^{-2\alpha_0z(\theta_t\w)}\(f(t,\me^{\alpha_nz(\theta_t\w)}v_{\al_n})
-f(t,\me^{\alpha_0z(\theta_t\w)}v_{\al_0}),
\me^{\al_nz(\theta_t\w)}v_{\al_n}-\me^{\al_0z(\theta_t\w)}v_{\al_0}\)\notag\\
:=&J_{31}+\me^{-2\alpha_0z(\theta_t\w)}(J_{32}+J_{33}).
\end{align}
Under the assumption \textbf{(F3)} and Young's inequality, we have for each $n\in\N$,
\begin{align}&\sum_{i\in\Z}|f_i(t,\me^{\alpha_nz(\theta_t\w)}v_{\al_n,i})|\notag\\
\leqslant&
\me^{\alpha_n(q-1)z(\theta_t\w)}\sum_{i\in\Z}|\psi_{3i}(t)||v_{\al_n,i}|^{q-1}
+\sum_{i\in\Z}|\psi_{4i}(t)|\notag\\
\leqslant&\me^{(q-1)\cZ}\sum_{i\in\Z}\(\frac{q-1}{q}|v_{\al_n,i}|^q
+\frac1q|\psi_{3i}(t)|^q\)+\|\psi_{4}(t)\|_1\notag\\
\leqslant&C\me^{(q-1)\cZ}\(\|v_{\al_n}\|_q^q+\|\psi_{3}(t)\|_q^q\)
+\|\psi_{4}(t)\|_1<+\8,\label{3.22}
\end{align}
which means that $f(t,\me^{\alpha_nz(\theta_t\w)}v_{\al_n})\in\ell^1$
for all $n\in\N$ and hence by \eqref{2.2},
it yields that
\begin{align}
J_{31}\leqslant&\cE_{-1,n}\|f(t,\me^{\alpha_nz(\theta_t\w)}v_{\al_n})\|_1
\|v_{\al_n}-v_{\al_0}\|\notag\\
\leqslant&2\cG_\de^{\frac12}\|f(t,\me^{\alpha_nz(\theta_t\w)}v_{\al_n})\|_1\cE_{-1,n}
\end{align}
and
\begin{align}
J_{32}\leqslant&\(\|f(t,\me^{\alpha_nz(\theta_t\w)}v_{\al_n})\|_1
+\|f(t,\me^{\alpha_0z(\theta_t\w)}v_{\al_0})\|_1\)
\cE_{1,n}\|v_{\al_n}\|\notag\\
\leqslant&\cG_\de^{\frac12}\(\|f(t,\me^{\alpha_nz(\theta_t\w)}v_{\al_n})\|_1
+\|f(t,\me^{\alpha_0z(\theta_t\w)}v_{\al_0})\|_1\)
\cE_{1,n}.
\end{align}
For $J_{33}$, we can use \textbf{(F1)} and obtain that
\begin{align}
J_{33}\leqslant&\|\psi_{2}(t)\|_\8\|\me^{\al_nz(\theta_t\w)}v_{\al_n}
-\me^{\al_0z(\theta_t\w)}v_{\al_0}\|^2\notag\\
\leqslant&\|\psi_{2}(t)\|_\8\(\(\me^{\al_nz(\theta_t\w)}-\me^{\al_0z(\theta_t\w)}\)\|v_{\al_n}\|
+\me^{\al_0z(\theta_t\w)}\|v_{\al_n}-v_{\al_0}\|\)^2\notag\\
\leqslant&2\|\psi_2(t)\|_{\8}\(\cG_\de\cE_{1,n}^2
+\me^{2\al_0z(\theta_t\w)}\|v_{\al_n}-v_{\al_0}\|^2\).\label{3.25}
\end{align}

Thus by setting $L_1$ and $L_2(n,t)$ such that
$$\frac12L_1:=C\nu_0\me^{(p-2)\cZ}\cG_{\de}^{\frac{p-2}{2}}+|\cZ-\lam|
+2\esssup_{t\in[\tau,T]}\|\psi_2(t)\|_{\8},$$
which is independent of $t$ and $n$ and
\begin{align}\frac12L_2(n,t):=&C\nu_0\cG_\de^{\frac{p}2}\cE_{p-2,n}+C\cZ\cG_\de|\al_n-\al_0|
+2\me^{2\cZ}\cG_\de\esssup_{t\in[\tau,T]}\|\psi_2(t)\|_{\8}\cE_{1,n}^2\notag\\
&+2\cG_\de^{\frac12}\(C\me^{(q-1)\cZ}\(\|v_{\al_n}\|_q^q+\|\psi_{3}(t)\|_q^q\)+\|\psi_{4}(t)\|_1\)
\(\cE_{-1,n}+\me^{2\cZ}\cE_{1,n}\),\label{3.26}
\end{align}
one can infer from \eqref{3.16} to \eqref{3.25} that
\be\label{3.27}\frac{\di}{\di t}\|v_{\al_n}-v_{\al_0}\|^2
\leqslant L_1\|v_{\al_n}-v_{\al_0}\|^2+L_2(n,t).\ee
Applying Gronwall's inequality to \eqref{3.27}, we can obtain
\begin{align}\label{3.28}\|v_{\al_n}(t)-v_{\al_0}(t)\|^2\leqslant&
\|v_{\al_n}(\tau)-v_{\al_0}(\tau)\|^2\me^{L_1(t-\tau)}
+\int_\tau^tL_2(n,s)\me^{L_1(t-s)}\di s\notag\\
\leqslant&\left|\me^{-\al_nz(\theta_\tau\w)}-\me^{-\al_0z(\theta_\tau\w)}\right|^2
\|u_\tau\|^2\me^{L_1(T-\tau)}
+\int_\tau^TL_2(n,s)\me^{L_1(T-s)}\di s.\end{align}

Then combining \eqref{3.11} and \textbf{(F3)},
we can obtain the uniform boundedness of the second term (integral)
on the right side of \eqref{3.28} for all $n\in\N^+$.
Moreover, each term in \eqref{3.26} contains either $\cE_{r,n}$
for some $r\in\R$ or $|\al_n-\al_0|$,
which does not depend on the time variable $t$.
This implies that as $n\ra\8$,
$$\int_\tau^TL_2(n,s)\me^{L_1(t-s)}\di s\ra0,$$
and hence the right side of \eqref{3.28} tends to $0$ as $n\ra\8$.
The proof is thus finished.
\eo

Based on the uniform continuity of $\vp_\al$ with respect to $\al$ given in Lemma \ref{le3.7}
and the continuity of $\vp_\al$ with respect to the initial datum,
one can trivially show that as $n\ra\8$,
\begin{align*}&\|\vp_{\al_n}(t,\tau,\w,u_n)-\vp_{\al_0}(t,\tau,\w,u_0)\|\\
\leqslant&\|\vp_{\al_n}(t,\tau,\w,u_n)-\vp_{\al_0}(t,\tau,\w,u_n)\|
+\|\vp_{\al_0}(t,\tau,\w,u_n)-\vp_{\al_0}(t,\tau,\w,u_0)\|\ra0,\end{align*}
whenever $\al_n\ra\al_0$ and $u_n\ra u_0$ in $\ell^2$ as $n\ra\8$.
This asserts the condition (1) of Theorem \ref{th3.2}.

Next we check the condition (2) of Theorem \ref{th3.2},
for which we need the following conclusion.

\bt\label{th3.8} Let $ K_n $ be a sequence of subsets of $X$ and $K\subset X$.
Suppose that for each subsequence $n_k$ of $n$,
and $x_{n_k} \in K_{n_k}$ with $x_{n_k}\ra x_0\in X$,
it holds that $x_0\in K$. Then
\be\label{3.29}\Cap_{n\geqslant 1}\ol{\Cup_{N\geqslant n} K_N}^X\subset K.\ee
\et
\bo Take $x_0$ in the left side of \eqref{3.29}. Then for each $n\geqslant1$,
$$x_0\in\ol{\Cup_{N\geqslant n} K_N}^X.$$
This  indicates that for each $k\geqslant1$,
there exists $n_k\geqslant1$ and $x_{n_k}\in K_{n_k}$ such that
$$\|x_{n_k}-x_0\|<\frac{1}{k}.$$
Furthermore, $n_k$ can be chosen strictly increasing, i.e., $n_{k+1}>n_k$.
Obviously, $x_{n_k}\ra x_0$ and by the supposition, $x_0\in K$. The proof is finished.\eo

Recall the results in \cite{WJ24}.
We know that for each $\al\in\R$,
the process $\vp_{\al}$ has a pullback $\cD_\al$-absorbing set
$K_\al=\{K_\al(\tau,\w)\}_{(\tau,\w)\in\R\X\W}\in\cD_\al$ with
$$K_\al(\tau,\w)=\{u\in\ell^2:\|u\|^2\leqslant\me^{2\al z(\w)}\cR(\al,\tau,\w)\},$$
where
$$\cR(\al,\tau,\w):=1+2\int_{-\8}^0\me^{\lam_0s
+2\al\int_s^0z(\theta_r\w)\di r-2\al z(\theta_s\w)}
\|\psi_1(\tau+s)\|_1\di s.$$
\bl\label{le3.9}Let $(\tau,\w)\in\R\X\W$ be arbitrary,
$\{\al_n\}_{n\in\N}$ be a sequence with $\al_n\ra\al_0$ as $n\ra\8$
and $\{x_n\}_{n\in\N}$ be a sequence in $\ell^2$ with $x_n\in K_{\al_n}(\tau,\w)$ for $n\in\N^+$
such that $x_n\ra x_0\in\ell^2$. Then $x_0\in K_{\al_0}(\tau,\w)$.
\el
\bo Indeed, for each fixed $(\tau,\w)\in\R\X\W$, we can let $g:\ell^2\X\R\ra\R$ such that
$$g(x,\al)=\me^{2\al z(\w)}\cR(\al,\tau,\w)-\|x\|^2.$$
By the Dominated Convergence Theorem, it is easy to see that
$\cR(\al,\tau,\w)$ is continuous in $\al$,
and hence $g(x,\al))$ is joint continuous in $(x,\al)$.

Note that for $x_n$ and $\al_n$, we have
$$g(x_n,\al_n)\geqslant0.$$
Now let $n\ra\8$. Then by the locally sign-preserving property of continuous functions, we know
$$g(x_0,\al_0)\geqslant0,\hs\mb{i.e.,}\hs
\me^{2\al_0 z(\w)}\cR(\al,\tau,\w)-\|x_0\|^2\geqslant0.$$
This indicates that $x_0\in K_{\al_0}(\tau,\w)$ and ends the proof.
\eo

Theorem \ref{th3.8} and Lemma \ref{le3.9} have actually guaranteed the condition (2) of
Theorem \ref{th3.2}.
To check the condition (3),
we can verify the local uniform pullback $\cD_\al$-asymptotic compactness
of $\{\vp_\al\}_{\al\in\R}$ by Theorem \ref{th3.4A}.
For this, we need again to recall the estimations in \cite{WJ24} as follows.
\bl\label{le3.10} Let the assumptions \textbf{(F1)}, \textbf{(F2)} and \textbf{(F3)}
hold for $f$, $\al_0\in\R$
with a neighborhood $\cI$ in $\R$ and $\me^{\al z(\theta_{-t}\w)}v_{\tau-t}\in D_\al(\tau,\w)$
with $D_\al\in\cD_\al$ for each $\al\in\cI$.
Then for each $\tau\in\R$, $\bP$-a.s. $\w\in\W$, $\al\in\cI$ and $\ve>0$,
there exist $T=T(\ve,\tau,\w,\cI,\{D_\al\}_{\al\in\cI})>0$ and
$N=N(\ve,\tau,\w,\cI,\{D_\al\}_{\al\in\cI})\in\N$ such that
when $t>T$, it holds that
\be\label{3.30}\sum_{|i|>N}|v_{\al,i}(\tau,\tau-t,\theta_{-\tau}\w,v_{\al,\tau-t})|^2<\ve,\ee
\be\label{3.31}\mb{and}\hs\sum_{|i|>N}|v_{\al,i}
(\tau,\tau-t,\theta_{-\tau}\w,v_{\al,\tau-t})|^q<\ve.\ee
\el
\bo For the proofs of \eqref{3.30} and \eqref{3.31},
we can refer to Lemmas 3.7 and 3.10 of \cite{WJ24}.
The only differences between the conclusions \eqref{3.30}, \eqref{3.31}
and \cite[Lemmas 3.7 and 3.10]{WJ24}
lie in the dependence on the parameter $\al$: the former depends
on the neighborhood $\cI$ of $\al_0$ uniformly,
while the latter depends on $\al$ pointwise.

Following the proofs of Lemmas 3.7 and 3.10 in \cite{WJ24} in the framework of this article,
one can see that every $\pm z(\theta_s\w)$ (resp. $\pm z(\theta_\sig\w)$,
$\pm z(\theta_\vsig\w)$ etc.)
should appear with a coefficient $\al$. If we take
$$\hat\al:=\sup_{\al\in\cI}|\al|$$
and let $\hat\al|z(\theta_s\w)|$ (resp. $\hat\al|z(\theta_\sig\w)|$,
$\hat\al|z(\theta_\vsig\w)|$ etc.) replace
$\pm\al z(\theta_s\w)$ (resp.  $\pm\al z(\theta_\sig\w)$, $\pm\al z(\theta_\vsig\w)$ etc.),
the procedures of the proofs for Lemmas 3.7 and 3.10 in \cite{WJ24} can work smoothly with
a slightly modified usage of \eqref{2.5}.
As a result, these procedures can deduce the uniform conclusions \eqref{3.30} and \eqref{3.31}.
Since the argument is almost the same as Lemmas 3.7 and 3.10 in \cite{WJ24} but rather tedious,
we omit the details here.
\eo

The property given in Lemma \ref{le3.10} can be called
\emph{local uniform pullback $\cD_\al$-asymptotic nullity}
of $\{\vp_\al\}_{\al\in\R}$.
The local uniform pullback $\cD_\al$-asymptotic nullity of $\{\vp_\al\}_{\al\in\R}$ can similarly
imply the local uniform pullback $\cD_\al$-asymptotic compactness,
and hence we have proved the condition (3$'$) given in Theorem \ref{th3.4A}.

Summarizing all the discussions above, we have actually proved the first main theorem.
\bt Let the assumptions \textbf{(F1)}, \textbf{(F2)} and \textbf{(F3)} hold for $f$
and $\cA_\al$ be the pullback random $(\ell^2,\ell^q)$-attractor
with respect to $\cD_{\al}$ for $\vp_\al$
with $q\in[1,\8)$.
Then for every $(\tau,\w)\in\R\X\W$ and $\al_0\in\R$, it holds that
$$\lim_{\al\ra\al_0}{\rm dist}_{\ell^2\cap\ell^q}(\cA_{\al}(\tau,\w),\cA_{\al_0}(\tau,\w))=0.$$
\et

\section{Convergence of invariant sample measures}\label{s4}
\subsection{Preliminaries on invariant sample measures}

We first give the relative definitions and theorems of invariant sample measures
(with two variables) for NRDS's; see \cite{WZ23,CY23}.

Let $X$ be a Banach space with norm $\|\cdot\|$ and
$\vp$ be an NRDS over $(\W,\cF,\bP,\{\theta_t\}_{t\in\R})$.
A family of Borel probability measures $\{\mu_{\tau,\w}\}_{(\tau,\w)\in\R\X\W}$ on $X$
is called an \emph{invariant sample measures} for $\vp$,
if for each $t\geqslant0$, $\tau\in\R$, $\bP$-a.s. $\w\in\W$ and all $U\in\cB(X)$,
it holds that
\be\label{4.1}
\mu_{\tau+t,\theta_t\w}(U)=\mu_{\tau,\w}((\vp(t,\tau,\w,\cdot))^{-1}U).\ee

The construction of the invariant sample measures usually involves
the \emph{generalized Banach limit}, which is any linear functional,
denoted by $\disp\LIM_{t\ra-\8}$ (we only need to consider the case $t\ra-\8$,
see more in \cite{WZC20,WZ23,WWL22}), defined for an arbitrary bounded real-valued function
on $(-\8,a]$ for some $a\in\R$ and satisfying
\benu
\item[(1)] $\disp\LIM_{t\ra-\8}\zeta(t)\geqslant0$ for nonnegative functions $\zeta$
on $(-\8,a]$;
\item[(2)] $\disp\LIM_{t\ra-\8}\zeta(t)=\lim_{t\ra-\8}\zeta(t)$ if the latter limit exists.
\eenu

Let $(\fA,\rho)$ be a metric space.
For each $\al\in\fA$, let $\vp_\al$ be an NRDS over $(\W,\cF,\bP,\{\theta_t\}_{t\in\R})$
on the state space $X$ and $\cD_\al$ be a universe for $\vp_\al$.
We introduce \cite[Theorem 4.1]{CY23} to discuss the convergence property
of invariant sample measures for the family $\{\vp_\al\}_{\al\in\fA}$.

\bt\label{th4.1}
Suppose that for each $\al\in\fA$, $\vp_{\al}$ has a pullback random $\cD_\al$-attractor
$\cA_{\al}(\tau,\w)$ and a family of Borel probability measures
$\{\mu_{\tau,\w}^{\al}\}_{(\tau,\w)\in\R\X\W}$ such that $\mu_{\tau,\w}^{\al}$
is an invariant sample measures and is supported on $\cA_{\al}(\tau,\w)$.
Let $\al_0\in\fA$.
Suppose that
\benu\item[(1)]the union $\Cup_{\al\in\fU}\cA_{\al}(\tau,\w)$
is precompact in $X$ for $\tau\in\R$, $\bP$-a.s. $\w\in\W$ and
a certain neighborhood $\fU$ of $\al_0$ in $\fA$; and
\item[(2)]for each compact set $K\subseteq X$, $\tau\in\R$,
$\bP$-a.s. $\w\in\W$ and each sequence $\{\al_n\}_{n=1}^{\8}$ in $\fA$
with $\al_n\ra\al_0$ as $n\ra\8$,
the following equality holds:
\be\label{4.2A}
\lim_{n\ra+\8}\sup_{u\in K}
\|\vp_{\al_n}(t,\tau-t,\theta_{-t}\w,u)-\vp_{\al_0}(t,\tau-t,\theta_{-t}\w,u)\|=0.
\ee
\eenu
Then for each $\tau\in\R$ and $\bP$-a.s. $\w\in\W$, there exists a subsequence
$\{\al_{n_k}\}_{k=1}^{+\8}$ of $\{\al_n\}_{n=1}^{+\8}$ depending on $\tau$ and $\w$,
and an invariant sample measure $\mu_{\tau,\w}^{\al_0}$ of $\vp_{\al_0}$ such that
$\mu_{\tau,\w}^{\al_{n_k}} $ converges weakly to $\mu_{\tau,\w}^{\al_0}$ as $k\ra+\8$.
\et

\subsection{Convergence of invariant sample measures}

First we recall the existence theorem of invariant sample measures for NRDS's
in \cite{WJ24}.
To this end, we need another assumption on $f$:

\noindent\textbf{(F4)} The nonlinear forcing $f$ satisfies \eqref{3.1A} with
\be\label{4.pp}\psi_3,\,\psi_4\in L_{\rm loc}^{2}(\R,\ell^{2}).\ee

\bt\label{th4.2}
Suppose the assumption \textbf{(F1)}, \textbf{(F2)} and \textbf{(F4)} hold for $f$.
Let $\cA_\al$ be the unique $(\ell^2,\ell^2)$-pullback random $\cD_\al$-attractor of $\vp_\al$
given in Theorem \ref{th3.8}.
Then for a given generalized Banach limit $\disp\LIM_{t\ra-\8}$ and a mapping $\xi:\R\X\W\ra\ell^2$
with $\{\xi(\tau,\w)\}_{(\tau,\w)\in\R\X\W}\in\cD_\al$ such that
$\xi(\tau,\theta_\tau\w)$ is continuous in $\tau$,
there exists a family of Borel probability measures $\{\mu^{\al}_{\tau,\w}\}_{(\tau,\w)\in\R\X\W}$
on $\ell^2$ such that the support of the measure $\mu^{\al}_{\tau,\w}$ is
contained in $\cA_\al(\tau,\w)$
and for all $\Upsilon\in\cC(\ell^2)$,
\begin{align}&
\LIM_{\tau\ra-\8}\frac{1}{t-\tau}\int_{\tau}^{t}\Upsilon
(\vp_\al(t-s,s,\theta_s\w,\xi(s,\theta_s\w)))\di s
=\int_{\cA_\al(t,\theta_t\w)}\Upsilon(u)\di\mu^{\al}_{t,\theta_t\w}(u)\notag\\
=&
\int_{\ell^2}\Upsilon(u)\di \mu^{\al}_{t,\theta_t\w}(u)
=\LIM_{\tau\ra-\8}\frac{1}{t-\tau}\int_{\tau}^{t}\int_{\ell^2}\Upsilon(\vp_\al(t-s,s,\theta_s\w,u))
\di\mu^{\al}_{s,\theta_s\w}(u)\di s.
\label{4.15}\end{align}
Additionally, $\mu^{\al}_{t,\w}$ is invariant in the sense that for all $\tau\in\R$ and
$\bP$-a.s. $\w\in\W$,
\be\label{4.16}
\int_{\cA_\al(\tau+t,\theta_t\w)}\Upsilon(u)\di \mu^{\al}_{\tau+t,\theta_t\w}(u)=
\int_{\cA_\al(\tau,\w)}\Upsilon(\vp_\al(t,\tau,\w,u))\di\mu^{\al}_{\tau,\w}(u),
\hs\mb{for all }t\geqslant0.
\ee\et

By using the definition \eqref{2.10} for $\vp_\al$, one can translate \eqref{4.16} into the version
for the solution $u_\al$ of \eqref{2.4} in the following form,
\be\label{4.17}
\int_{\cA_\al(\tau+t,\theta_t\w)}\Upsilon(u)\di \mu^\al_{\tau+t,\theta_t\w}(u)=
\int_{\cA_\al(\tau,\w)}\Upsilon(u_\al(t+\tau,\tau,\theta_{-\tau}\w,u))\di\mu^{\al}_{\tau,\w}(u),
\hs\mb{for all }t\geqslant0.
\ee

We aim to apply Theorem \ref{th4.1} to prove the convergence property of
the invariant sample measures given by Theorem \ref{th4.2}.
By the discussion above, we actually only need to check the condition (2)
in Theorem \ref{th4.1}.
For this purpose, we recall the continuity given by Lemma \ref{le3.7}.
Combining the definition \eqref{2.10}, we only need to ensure the following conclusion.

\bl\label{le4.3} Let the assumptions \textbf{(F1)}, \textbf{(F2)} and
\textbf{(F3)} hold for $f$.
Then for each compact set $K\subset\ell^2$, $\tau\in\R$,
$\bP$-a.s. $\w\in\W$ and $\al_n\in \fA$ with $\al_n\ra\al_0$ as $n\ra\8$,
the following equality holds true:
\be\label{4.3A}
\lim_{n\ra+\8}\sup_{u\in K}\|u_{\al_n}(t,\tau-t,\theta_{-\tau}\w,u)
-u_{\al_0}(t,\tau-t,\theta_{-\tau}\w,u)\|=0,
\ee
for all $t\geqslant0$.
\el

\bo We first fix $u\in\ell^2$, $\tau\in\R$, $t\geqslant0$ and $\bP$-a.s. $\w\in\W$.
By Lemma \ref{le3.7}, we know for each $\ve>0$, there exists an open neighborhood $\cO_{u}$
of $u$ in $\ell^2$ and $N_u\in\N^+$ such that for all $\phi\in\cO_u$ and $n>N_u$,
\be\label{4.4A}\|u_{\al_n}(t,\tau-t,\theta_{-\tau}\w,\phi)
-u_{\al_0}(t,\tau-t,\theta_{-\tau}\w,\phi)\|<\ve.\ee

For the compact set $K$ in $\ell^2$, the family $\{\cO_u\}_{u\in K}$ constitutes an open
covering of $K$ and hence has a finite sub-covering $\{\cO_{u^{(i)}}\}_{i=1}^k$ with
$u^{(i)}\in K$ for each $i\in[1,k]$ such that
$$K\subset\Cup_{i=1}^k\cO_{u^{(i)}}.$$
Let
$$N_K=\max\{N_{u^{(1)}},\,N_{u^{(2)}},\,\cdots,\,N_{u^{(k)}}\}.$$
Then one can easily see that when $n>N_K$,
\eqref{4.4A} is valid for each $\phi\in K$,
which obviously indicates \eqref{4.3A}.
The proof is complete.
\eo

Now we observe that the condition (1) in Theorem \ref{th4.1} follows from
Lemma \ref{le3.10} and Theorem \ref{th3.4A}.
Furthermore, Lemma \ref{le4.3} and \eqref{2.10} immediately guarantee the condition
(2) in Theorem \ref{th4.1}.
Consequently, we can thereby obtain the main theorem in this section.

\bt\label{th4.4} Suppose the assumptions \textbf{(F1)}, \textbf{(F2)}, \textbf{(F3)} and
\textbf{(F4)} hold for $f$.
Let $\vp_\al$ be the NRDS defined as \ref{2.10} generated by \eqref{2.4}
and $\{\mu^{\al}_{\tau,\w}\}_{(\tau,\w)\in\R\X\W}$ be the invariant sample measures
of $\vp_\al$.
Then for each $\tau\in\R$, $\bP$-a.s. $\w\in\W$
and each sequence $\{\al_n\}_{n\in\N^+}$ in $\R$
with $\al_n\ra\al_0$ as $n\ra\8$,
there exists a subsequence $\{\al_{n_k}\}_{k\in\N^+}$ of $\{\al_n\}_{n\in\N^+}$
depending on $\tau$ and $\w$,
and an invariant sample measure $\mu_{\tau,\w}^{\al_0}$ of $\vp_{\al_0}$ such that
$\mu_{\tau,\w}^{\al_{n_k}}$ (or $\{\mu_{t,\theta_{t-\tau}\w}^{\al_{n_k}}\}_{t\in\R}$)
converges weakly to $\mu_{\tau,\w}^{\al_0}$ (or $\{\mu_{t,\theta_{t-\tau}\w}^{\al_{0}}\}_{t\in\R}$)
as $k\ra+\8$ in the sense that
$$\mu_{t,\theta_{t-\tau}\w}^{\al_{n_k}}\hs\mb{converges weakly to}\hs\mu_{t,\theta_{t-\tau}\w}^{\al_0}
\hs\mb{for every }t\in\R.$$
\et

\section{Convergence of stochastic Liouville type equation}\label{s5}
\subsection{Stochastic Liouville type equation}\label{ss5.1}
In this subsection, we first check that the invariant sample measures
$\{\mu_{\tau,\w}\}_{(\tau,\w)\in\R\X\W}$
obtained in Theorem \ref{th4.2} for $\vp$ fulfill
a stochastic Liouville type theorem (see \cite{CY23,WZL24}).

We adopt the stochastic Liouville type equation (see \cite{CY23,WZL24}) with respect to
the invariant sample measures $\{\mu_{\tau,\w}\}_{(\tau,\w)\in\R\X\W}$.
Let
$$V=\ell^2\cap\ell^q,\hs H=\ell^2\hs\mb{and}\hs V^*=\ell^2+\ell^{q'},$$
where $q'$ is the dual index of $q$ with $q'=\frac{q}{q-1}$ for $q=1$ and $q'=\8$ for $q=1$.
We set the norm of $V$ as $\|u\|_V=\|u\|+\|u\|_q$.
Denote the dual pairing between $v^*\in V^*$ and $v\in V$ by $\langle v^*,v\rangle$.
Let
\be\label{5.0A}\tilde{F}(t,u)=-\nu(t)Au+f(t,u)-\lam u.\ee
We know $\tilde{F}:\R\X V\ra V^*$.
Thus the equation \eqref{2.4} can be written as
\be\label{5.1}\di u=\tilde{F}(t,u)\di t+\al u\circ\di W.\ee

Let $\cT$ denote the class of test functions $\Psi:H\ra\R$ for \eqref{5.1}
to be bounded on bounded subsets of $H$ such that the following conditions hold:
\benu\item[(1)] for each $u\in V$, the Fr\'echet derivative $\Psi'(u)$ taken
in $H$ along $V$ exists.
More precisely, for each $u\in V$, there exists an element in $H$ denoted by $\Psi'(u)$ such that
$$\frac{|\Psi(u+v)-\Psi(v)-(\Psi'(u),v)|}{\|v\|}\ra0\hs\mb{as }\|v\|\ra0,\;v\in V;$$
\item[(2)] for each $u\in V$, $\Psi'(u)\in V$ and the mapping $u\mapsto\Psi'(u)$ is continuous
and bounded as a functional from $V$ into $V$;
\item[(3)] for each $u\in V$, the second-order Fr\'echet derivative $\Psi''(u)$ is
a bounded bilinear operator from $H\X H$ to $\R$ and the mapping $u\mapsto\Psi''(u)$ is continuous
and bounded from $V$ into $L(H\X H,\R)$;
\item[(4)] for every global solution $u(t)$ of \eqref{5.1}, the following It\^o's Formula holds for
all $t\geqslant s$ and $\bP$-a.s. $\w\in\W$,
\begin{align}\Psi(u(t))-\Psi(u(s))=&\int_s^t\langle\tilde{F}(\vsig,u(\vsig)),
\Psi'(u(\vsig))\rangle\di\vsig
+\al\int_s^t\(u(\vsig),\Psi'(u(\vsig))\)\di W(\vsig)\notag\\
&+\frac{\al^2}2\int_s^t\Psi''(u(\vsig))\(u(\vsig),u(\vsig)\)\di\vsig.
\label{2.8n}\end{align}
\eenu
The class $\cT$ is not empty for \eqref{5.1}, and one can see examples in \cite{CY23,WZL24}.

For the following discussion, we set one more assumption for $f(t,u)$ as follows,

\noindent\textbf{(F5)} there is $\psi_5(t)=(\psi_{5i}(t))_{i\in\Z}$, such that
\be\label{3.58z}\left|\frac{\pa f_i}{\pa u_i}(t,u_i)\right|
\leqslant\psi_{5i}(t)\mb{ with }\psi_5\in L_{\rm loc}^\infty(\R;\ell^\infty).\ee

To prove the final conclusion, we need the following results.

\bl\label{le5.1} Suppose that the assumptions \textbf{(F1)}, \textbf{(F2)} and \textbf{(F3)}
hold for $f$.
Let $\cA_\al$ be the pullback random $(\ell^2,\ell^q)$-attractor of $\vp_\al$
with respect to $\cD_\al$ given in Theorem \ref{th3.5}.
Then for each $\al\in\R$, $\tau\in\R$, $\bP$-a.s. $\w\in\W$,  and $s,t\in\R$ with $s\leqslant t$,
the set
\be\label{5.1a}\cU_\al[s,t]:=\Cup_{\sig\in[s,t]}\cA_{\al}(\sig,\theta_{\sig-\tau}\w)\ee
is compact in $\ell^2$.\el
\bo We pick a sequence $\{\phi_k\}_{k\in\N^+}$ in $\cU_\al[s,t]$
arbitrarily with $\phi_k\in\cA_{\al}(\sig_k,\theta_{\sig_k-\tau}\w)$,
and it suffices to show that this sequence has a convergent subsequence
with its limit in $\cU_\al[s,t]$.

Since $\sig_k\in[s,t]$, we can assume that $\sig_k\ra\sig_0\in[s,t]$ as $k\ra\8$
(selecting its subsequence if necessary).
Take
$$\sig_*=\inf_{k\in\N}\sig_k\in[s,t].$$
Then by the invariance of the attractor, we have
$\tilde{\phi}_k\in\cA_{\al}(\sig_*,\theta_{\sig_*-\tau}\w)$ for each $n\in\N$ such that
$$\phi_k=u_{\al}(\sig_k,\sig_*,\theta_{-\tau}\w,\tilde{\phi}_k).$$
Moreover, since $\tilde{\phi}_{k}\in\cA_{\al}(\sig_*,\theta_{\sig_*-\tau}\w)$,
which is compact in $\ell^2$,
we can furthermore assume that
$$\tilde{\phi}_k\ra\phi_*\in\cA_\al(\sig_*,\theta_{\sig_*-\tau}\w),\hs
\mb{as }n\ra\8.$$
Recalling the condition (4) of a continuous NRDS in Section \ref{s2}
(or \cite[Theorem 2.5]{WJ24}),
we can infer from $\sig_k\ra\sig_0$ and $\tilde{\phi}_k\ra\phi_*$ that
$$
\|\phi_k-u_{\al}(\sig_0,\sig_*,\theta_{-\tau}\w,\phi_*)\|\leqslant
\|u_{\al}(\sig_k,\sig_*,\theta_{-\tau}\w,\tilde{\phi}_k)
-u_{\al}(\sig_0,\sig_*,\theta_{-\tau}\w,\phi_*)\|\ra0,$$
as $n\ra\8$.
The compactness is thus proved.
\eo

Let $U\subset H$.
We say that a mapping $U\ni u\mapsto\tilde\Psi(u)\in X$ for a Banach space $X$
is \textit{continuous in (norm of) $H$},
if for a sequence $\{u^{(n)}\}_{n\in\N^+}$
and $u$ in $U$,
\be\label{3.59z}
\|\tilde\Psi(u^{(n)})-\tilde\Psi(u)\|_X\ra0,\hs\mb{as }\|u^{(n)}-u\|\ra0.
\ee
\bl\label{le5.2} Suppose that the assumptions \textbf{(F1)}, \textbf{(F2)}, \textbf{(F3)},
\textbf{(F4)} and \textbf{(F5)} hold for $f$.
Let $\Psi\in\cT$ and $u\in\cA_\al(\tau,\w)$,
where $\al\in\R$ and $\cA_\al=\{\cA_\al(\tau,\w)\}_{(\tau,\w)\in\R\X\W}$ is the pullback random
$(\ell^2,\ell^q)$-attractor given by Theorem \ref{th3.5}.
Then the mappings
\be\label{5.aqa}\cA_{\al}(t,\theta_{t-\tau}\w)\ni u\mapsto\langle\tilde{F}(t,u),\Psi'(u)\rangle\ee
\be\label{5.awa}\cA_{\al}(t,\theta_{t-\tau}\w)\ni u\mapsto\Psi''(u)(u,u)
\hs\mb{and}\hs\cA_{\al}(t,\theta_{t-\tau}\w)\ni u\mapsto(u,\Psi'(u))\ee
are continuous in (norm of) $H$.
Moreover, for the invariant sample measures $\{\mu_{\tau,\w}\}_{(\tau,\w)\in\R\X\W}$ given
in Theorem \ref{th4.2}, the mappings
\be\label{3.60z}
\vsig\mapsto\int_{H}\langle\tilde{F}(\vsig,u),\Psi'(u)\rangle
\mu_{\vsig,\theta_{\vsig-\tau}\w}(\di u)\hs\mb{and}
\hs\vsig\mapsto\int_{H}\Psi''(u)(u,u)\mu_{\vsig,\theta_{\vsig-\tau}\w}(\di u)
\ee
belong to $L^1_{\rm loc}(\R)$;
and the mapping
\be\label{5.AAA}\vsig\mapsto\int_{H}(u,\Psi'(u))\mu_{\vsig,\theta_{\vsig-\tau}\w}(\di u)\ee
belongs to $L^2_{\rm loc}(\R)$.
\el

\bo We first prove the results for $\langle\tilde{F}(t,u),\Psi'(u)\rangle$.
To this end, we first consider the continuity of the mapping
\be\label{3.61z}V\ni u\mapsto\tilde{F}(t,u)\in V^*\ee
in norm of $H$ for each $t\in\R$.
Take $u$, $v$, $w\in V$ with $\|u\|$, $\|v\|\leqslant R$ for some $R>0$
and set $\tilde{w}=u-v$.
For each fixed $t\in\R$, we observe by \eqref{3.58z} and
the Differential Mean-Value Theorem that
\begin{align*}&\left|\langle\tilde{F}(t,u)-\tilde{F}(t,v),w\rangle\right|\\
\leqslant&\nu_0\left|\langle|Bu|^{p-2}\otimes(Bu)-|Bv|^{p-2}\otimes(Bv),
Bw\rangle\right|+\lam\|\tilde{w}\|\|w\|+\sum_{i\in\Z}|f_i(t,u_i)-f_i(t,v_i)||w_i|\\
\leqslant&2(p-1)(\|Bu\|+\|Bv\|)^{p-2}\|B\tilde{w}\|\|w\|+\lam\|\tilde{w}\|\|w\|
+\(\sum_{i\in\Z}\left|\frac{\pa f_i}{\pa u_i}(t,\xi_i)\right|^2
|\tilde{w}_i|^2\)^{\frac12}\|w\|\\
\leqslant&\left[(p-1)4^{p-1}R^{p-2}+\lam+\|\psi_5(t)\|_\8\right]\|\tilde{w}\|\|w\|,
\end{align*}
where we also used the trivial estimate (14) of \cite{WJ24}.
Then we know that
\begin{align*}\|\tilde{F}(t,u)-\tilde{F}(t,v)\|_{V^*}
=&\sup_{\|w\|_V\ne0}\frac{\left|\langle\tilde{F}(t,u)-\tilde{F}(t,v),
w\rangle\right|}{\|w\|_V}\\
\leqslant&\sup_{\|w\|_V\ne0}\frac{\left|\langle\tilde{F}(t,u)-\tilde{F}(t,v),
w\rangle\right|}{\|w\|}\\
\leqslant&\left[(p-1)4^{p-1}R^{p-2}+\lam+\|\psi_5(t)\|_\8\right]\|\tilde{w}\|,
\end{align*}
which indicates the continuity \eqref{3.61z} in $H$ (also $V$ of course) immediately.

Then we consider the continuity of the mapping $u\mapsto\Psi'(u)$
from $\cA_\al(t,\theta_{t-\tau}\w)$ to $V$ in norm of $H$.
Let $\{u^{(n)}\}_{n\in\N}$ be a sequence in $\cA_\al(t,\theta_{t-\tau}\w)$ and
$\|u^{(n)}-u\|\ra0$ with $u\in H$ as $n\ra\8$.
Noting that $\cA_\al(t,\theta_{t-\tau}\w)$ is compact in the spaces $\ell^2$ and $\ell^q$, we know that
$u\in\cA_\al(t,\theta_{t-\tau}\w)$.
And each subsequence $u^{(n_k)}$ of $u^{(n)}$ has a new subsequence
converging to $u^*$ in $\cA_\al(t,\theta_{t-\tau}\w)$ in norm (or topology) of $\ell^q$.
Since $(\ell^2,\ell^q)$ is limit-identical, we know $u^*=u$.
It is easy to see that $u^{(n)}\ra u$ in norm of $\ell^q$ and hence of $V$.
As a result, the condition (2) for the class $\cT$ implies
that
$$\|\Psi'(u^{(n)})-\Psi'(u)\|_V\ra0,$$
and therefore the continuity of $u\mapsto\Psi'(u)$ from $\cA_\al(t,\theta_{t-\tau}\w)$ to $V$
(also $H$ obviously) in norm of $H$.
Finally the continuity of $u\mapsto\langle\tilde{F}(t,u),\Psi'(u)\rangle$ follows immediately.
\vs

Then we show the mapping
$$\vsig\mapsto\int_{H}\langle\tilde{F}(\vsig,u),\Psi'(u)\rangle
\mu_{\vsig,\theta_{\vsig-\tau}\w}(\di u)$$
belongs to $L^1_{\rm loc}(\R)$.
Due to again the fact that for each $(\tau,\w)\in\R\X\W$,
$\cA_\al(\tau,\w)$ is compact in $H$ and that the invariant sample measures $\mu_{\tau,\w}$
are supported by $\cA_\al(\tau,\w)$,
the mapping \eqref{3.60z} is surely well defined.
Note also that
\begin{align}\left|\langle \tilde{F}(\vsig,u),\Psi'(u)\rangle\right|
\leqslant&\nu_0\left|\langle Au,\Psi'(u)\rangle\right|
+\lam\left|\( u,\Psi'(u)\)\right|
+\left|\langle f(\vsig,u),\Psi'(u)\rangle\right|\notag\\
\leqslant&2^p\nu_0\|u\|^{p-1}\|\Psi'(u)\|+\lam\|u\|\|\Psi'(u)\|\notag\\
&+\sum_{i\in\Z}\(\psi_{3i}(\vsig)|u_i|^{q-1}+\psi_{4i}(\vsig)\)|\Psi'_i(u)|\notag\\
\leqslant&\(2^p\nu_0\|u\|^{p-1}+\lam\|u\|+\|\psi_4(\vsig)\|_1\)\|\Psi'(u)\|_V\notag\\
&+\sum_{i\in\Z}\psi_{3i}(\vsig)\(\frac{q-1}{q}|u_i|^{q}+\frac1q|\Psi'(u)|^q\)\notag\\
\leqslant&\(2^p\nu_0\|u\|^{p-1}+\lam\|u\|+\|\psi_4(\vsig)\|_1\)\|\Psi'(u)\|_V\notag\\
&+\frac{q-1}{q}\|\psi_{3}(\vsig)\|_1\|u\|^{q}+\frac1q\|\psi_3(\vsig)\|_1\|\Psi'(u)\|_V^q.
\label{5.2a}\end{align}
Now we consider an arbitrary interval $[s,t]\subset\R$.
By (2) of the definition of $\cT$, we know $\|\Psi'(u)\|_V$ is bounded;
by the compactness of $\cU_\al[s,t]$, we also know that $\|u\|$ is bounded
on $\cU_\al[s,t]$.
Hence we have a positive constant $\sM_\al$ independent of $\vsig$ such that
$$\left|\langle \tilde{F}(\vsig,u),\Psi'(u)\rangle\right|
\leqslant\sM_\al(1+\|\psi_3(\vsig)\|_1+\|\psi_4(\vsig)\|_1)$$
and then
\begin{align}&\left|\int_{H}\langle \tilde{F}(\vsig,u),\Psi'(u)\rangle
\mu^{\al}_{\vsig,\theta_{\vsig-t}\w}(\di u)\right|
=\left|\int_{\cA_{\al_n}(\vsig,\theta_{\vsig-t}\w)}
\langle \tilde{F}(\vsig,u),\Psi'(u)\rangle
\mu^{\al}_{\vsig,\theta_{\vsig-t}\w}(\di u)\right|\notag\\
\leqslant&\sM_\al(1+\|\psi_3(\vsig)\|_1+\|\psi_4(\vsig)\|_1),
\label{5.aea}\end{align}
which ensures that the first mapping in \eqref{3.60z} belongs to $L^1_{\rm loc}(\R)$
by the assumption \textbf{(F3)}.
\Vs

The continuity of the remaining two mappings in \eqref{5.awa} is obvious,
as long as we notice that $\Phi'(u)$ is continuous from $\cA_\al(t,\theta_{t-\tau}\w)$
to $H$ and so is $\Phi''(u)$ from $\cA_\al(t,\theta_{t-\tau}\w)$
to $L(H\X H,\R)$ by the definition of the class $\cT$.

Next, we recall the definition of $\cT$ that $\Psi'(u)$ is bounded from $V$ to $V$
and $\Psi''(u)$ is bounded from $V$ to $L(H\X H,\R)$.
By the embedding $V\subset H$, we also know that $\Psi'(u)$ is bounded in $H$.
Combining the boundedness of $\|u\|$ on the compact set $\cU_\al[s,t]$, one can directly
see that $|(u,\Psi'(u))|$ and $|\Psi''(u)(u,u)|$ are uniformly bounded on $\cU_\al[s,t]$.
Then similar to the analysis in \eqref{5.aea}, we can obtain that the second mapping
in \eqref{3.60z} belongs to $L^1_{\rm loc}(\R)$ and the mapping \eqref{5.AAA}
belongs to $L^2_{\rm loc}(\R)$.
The proof is hence accomplished.\eo

\br Here the requirement that the integrand of the stochastic integral to
belong to $L^2_{\rm loc}(\R)$
is to guarantee that the stochastic integral makes sense; see \cite[Section 2.3]{K05}.\er

Now we are well prepared to show the stochastic Liouville type theorem.
\bt\label{th5.4} Suppose that the assumptions \textbf{(F1)}, \textbf{(F2)},
\textbf{(F3)}, \textbf{(F4)}
and \textbf{(F5)} hold for $f$.
Let $\cA_\al$ be the pullback random $(\ell^2,\ell^q)$-attractor of $\vp_\al$
with respect to $\cD_\al$ given in Theorem \ref{th3.5} and
$\{\mu^{\al}_{\tau,\w}\}_{(\tau,\w)\in\R\X\W}$ be the invariant sample measures constructed in
Theorem \ref{th4.2}.
Then the invariant sample measures $\{\mu^{\al}_{\tau,\w}\}_{(\tau,\w)\in\R\X\W}$ satisfy
the stochastic Liouville type equation \eqref{1.4}, for $s,t\in\R$ with $s\leqslant t$.\et
\bo Based on the definition of the class $\cT$,
and by \eqref{4.15}, \eqref{4.17} and \eqref{2.8n}, we know that for all $t\geqslant\tau$,
\begin{align}&\int_{H}\Psi(\phi)\mu^{\al}_{t,\theta_{t-\tau}\w}(\di\phi)
-\int_{H}\Psi(\phi)\mu^{\al}_{s,\theta_{s-\tau}\w}(\di\phi)\notag\\
=&\int_{\cA_\al(t,\theta_{t-\tau}\w)}\Psi(\phi)\mu^{\al}_{t,\theta_{t-\tau}\w}(\di\phi)
-\int_{\cA_\al(s,\theta_{s-\tau}\w)}\Psi(\phi)\mu^{\al}_{s,\theta_{s-\tau}\w}(\di\phi)\notag\\
=&\int_{\cA_\al(s,\theta_{s-\tau}\w)}\(\Psi(u_\al(t,s,\theta_{-\tau}\w,\phi))-\Psi(\phi)\)
\mu^{\al}_{s,\theta_{s-\tau}\w}(\di\phi)\notag\\
=&\int_{\cA_\al(s,\theta_{s-\tau}\w)}\int_s^t\langle\tilde{F}(\vsig,u_\al(\vsig)),
\Psi'(u_\al(\vsig))\rangle\di\vsig\cdot\mu^{\al}_{s,\theta_{s-\tau}\w}(\di\phi)\notag\\
&+\al\int_{\cA_\al(s,\theta_{s-\tau}\w)}\int_s^t
\(u_\al(\vsig),\Psi'(u_\al(\vsig))\)\di\tilde{W}(\vsig)\cdot
\mu^{\al}_{s,\theta_{s-\tau}\w}(\di\phi)\notag\\
&+\frac{\al^2}2\int_{\cA_\al(s,\theta_{s-\tau}\w)}\int_s^t\Psi''(u_\al(\vsig))
\(u_\al(\vsig),u_\al(\vsig)\)\di\vsig\cdot
\mu^{\al}_{s,\theta_{s-\tau}\w}(\di\phi)\notag\\
:=&M_1+M_2+M_3,\label{3.62z}
\end{align}
where $\tilde{W}(\vsig)=W(-\tau+\vsig)-W(-\tau)$ comes from
the stochastic datum $\theta_{-\tau}\w$.
Note by (3) of the definition of NRDS's and \eqref{2.10} that
\be\label{3.63z}
u_\al(\vsig,s,\theta_{-\tau}\w,u_\al(s,\sig,\theta_{-\tau}\w,\phi))
=u_\al(\vsig,\sig,\theta_{-\tau}\w,\phi).\ee
Then using \eqref{4.15}, \eqref{3.63z}, Lemma \ref{le5.2}, Fubini's Theorem and
the invariance of $\mu^{\al}_{\tau,\w}$ to $M_1$, we have
\begin{align}M_1=&\LIM_{\tau'\ra-\8}\frac{1}{s-\tau'}\int_{\tau'}^{s}\int_{H}\int_s^t
\langle\tilde{F}(\vsig,u_\al(\vsig,\sig,\theta_{-\tau}\w,\phi)),
\Psi'(u_\al(\vsig,\sig,\theta_{-\tau}\w,\phi))\rangle
\di\vsig\cdot \mu^{\al}_{\sig,\theta_{\sig-\tau}\w}(\di\phi)\di\sig\notag\\
=&\LIM_{\tau'\ra-\8}\frac{1}{s-\tau'}\int_{\tau'}^s\int_s^t\int_{H}
\langle\tilde{F}(\vsig,u_\al(\vsig,\sig,\theta_{-\tau}\w,\phi)),
\Psi'(u_\al(\vsig,\sig,\theta_{-\tau}\w,\phi))\rangle
\mu^{\al}_{\sig,\theta_{\sig-\tau}\w}(\di\phi)\di\vsig\di\sig\notag\\
=&\LIM_{\tau'\ra-\8}\frac{1}{s-\tau'}\int_{\tau'}^s\int_s^t\int_{H}
\langle\tilde{F}(\vsig,u),\Psi'(u)\rangle
\mu^{\al}_{\vsig,\theta_{\vsig-\tau}\w}(\di u)\di\vsig\di\sig\notag\\
=&\int_s^t\int_{H}\langle\tilde{F}(\vsig,u),\Psi'(u)\rangle
\mu^{\al}_{\vsig,\theta_{\vsig-\tau}\w}(\di u)\di\vsig,\label{3.64z}
\end{align}
where we have replaced $t$, $\tau$ and $\w$ in \eqref{4.17} by $\vsig-\sig$, $\sig$
and $\theta_{\sig-\tau}\w$, respectively, to obtain the third equality.
Similarly, for $M_2$ and $M_3$, by Lemma \ref{le5.2}, we also have
\begin{align}M_2=&\al\LIM_{\tau'\ra-\8}\frac{1}{s-\tau'}\suo\int_{\tau'}^s\suo\int_{H}\!\!\int_s^t
\(u_\al(\vsig,\sig,\theta_{-\tau}\w,\phi),\Psi'(u_\al(\vsig,\sig,\theta_{-\tau}\w,\phi))\)
\di\tilde{W}(\vsig)\cdot \mu^{\al}_{\sig,\theta_{\sig-\tau}\w}(\di\phi)\di\sig\notag\\
=&\al\int_s^t\int_{H}\(u,\Psi'(u)\)
\mu^{\al}_{\vsig,\theta_{\vsig-\tau}\w}(\di u)\di\tilde{W}(\vsig)\hs\mb{and}
\end{align}
\be M_3=\frac{\al^2}{2}\int_s^t\int_{H}\Psi''(u)\(u,u\)
\mu^{\al}_{\vsig,\theta_{\vsig-\tau}\w}(\di u)\di \vsig.\label{3.66z}
\ee
Then the stochastic Liouville type equation \eqref{1.4} follows from
the equalities from \eqref{3.62z} to \eqref{3.66z}.
The proof is thus finished.
\eo
\subsection{Convergence of stochastic Liouville type equation}

Since the stochastic Liouville type equation relies on the coefficient $\al$ of
the stochastic term, it is also an interesting problem to study what will happen
when a sequence $\{\al_n\}_{n\in\N^+}\subset\R$ tends to some number $\al_0$ in $\R$.
Actually, in Theorem \ref{th4.4}, we have obtained the weak convergence of a subsequence
of $\{\mu^{\al_n}_{\tau,\w}\}_{(\tau,\w)\in\R\X\W}$ for a sequence $\{\al_n\}_{n\in\N^+}$ in $\R$
with $\al_n\ra\al_0\in\R$ as $n\ra\8$.
Following this result, one can further deduce the following theorem.

\bt\label{th5.5} Suppose that the assumptions \textbf{(F1)}, \textbf{(F2)},
\textbf{(F3)}, \textbf{(F4)} and \textbf{(F5)} hold for $f$.
Let $\cA_\al$ be the pullback random $(\ell^2,\ell^q)$-attractor of $\vp_\al$
with respect to $\cD_\al$ given in Theorem \ref{th3.5} and
$\{\mu^{\al}_{\tau,\w}\}_{(\tau,\w)\in\R\X\W}$ be the invariant sample measures constructed in
Theorem \ref{th4.2}.
Then for each $\tau\in\R$, $\bP$-a.s. $\w\in\W$, each convergent sequence $\{\al_n\}_{n\in\N^+}$
in $\R$ with $\al_n\ra\al_0\in\R$,
there exists a subsequence $\{\al_{n_k}\}_{k\in\N^+}$ of $\{\al_n\}_{n\in\N^+}$
depending on $\tau$ and $\w$,
and an invariant sample measure $\mu_{\tau,\w}^{\al_0}$ of $\vp_{\al_0}$ such that
for each $\Psi\in\cT$ and $s,t\in\R$ with $s\leqslant t$,
the stochastic Liouville type equation
\begin{align}&\int_{H}\Psi(u)\mu^{\al_{n_k}}_{t,\theta_{t-\tau}\w}(\di u)
-\int_{H}\Psi(u)\mu^{\al_{n_k}}_{s,\theta_{s-\tau}\w}(\di u)\notag\\
=&\int_s^t\int_H\langle\tilde{F}(\vsig,u),\Psi'(u)\rangle
\mu^{\al_{n_k}}_{\vsig,\theta_{\vsig-\tau}\w}(\di u)\di\vsig
+\al_{n_k}\int_s^t\int_H(u,\Psi'(u))\mu^{\al_{n_k}}_{\vsig,\theta_{\vsig-\tau}\w}(\di u)
\di\tilde{W}(\vsig)\notag\\
&+\frac{\al_{n_k}^2}2\int_s^t\int_H\Psi''(u)(u,u)
\mu^{\al_{n_k}}_{\vsig,\theta_{\vsig-\tau}\w}(\di u)\di\vsig
\label{5.11B}\end{align}
converges termwise to the following stochastic Liouville type equation
\begin{align}&\int_{H}\Psi(u)\mu^{\al_0}_{t,\theta_{t-\tau}\w}(\di u)
-\int_{H}\Psi(u)\mu^{\al_0}_{s,\theta_{s-\tau}\w}(\di u)\notag\\
=&\int_s^t\int_H\langle\tilde{F}(\vsig,u),\Psi'(u)\rangle
\mu^{\al_0}_{\vsig,\theta_{\vsig-\tau}\w}(\di u)\di\vsig
+\al_0\int_s^t\int_H(u,\Psi'(u))\mu^{\al_0}_{\vsig,\theta_{\vsig-\tau}\w}(\di u)
\di\tilde{W}(\vsig)\notag\\
&+\frac{\al_0^2}2\int_s^t\int_H\Psi''(u)(u,u)
\mu^{\al_0}_{\vsig,\theta_{\vsig-\tau}\w}(\di u)\di\vsig.
\label{5.12B}\end{align}
\et

To prove Theorem \ref{th5.5}, we still need some auxiliary results.

\bl\label{le5.6} Under the conditions of Theorem \ref{th5.5},
for each $\tau\in\R$, $\bP$-a.s. $\w\in\W$, $s,t\in\R$ with $s\leqslant t$
and a sequence $\{\al_n\}_{n\in\N^+}$ with $\al_n\ra\al_0\in\R$ as $n\ra\8$,
the set
$$\cU[s,t]:=\Cup_{n=0}^{\8}\cU_{\al_n}[s,t]$$
is precompact in $\ell^2$, where $\cU_{\al}[s,t]$ is defined in \eqref{5.1a}.\el
\bo We pick a sequence $\{\phi_k\}_{k\in\N^+}$ in $\cU[s,t]$
arbitrarily with $\phi_k\in\cA_{\al_{n_k}}(\sig_k,\theta_{\sig_k-\tau}\w)$,
and it suffices to show that this sequence has a convergent subsequence
with its limit in $\ell^2$.

Similar to the proof of Lemma \ref{le5.1}, we assume that $\sig_k\ra\sig_0\in[s,t]$ as $k\ra\8$
and take
$$\sig_*=\inf_{k\in\N}\sig_k\in[s,t].$$
Then we have
$\tilde{\phi}_k\in\cA_{\al_{n_k}}(\sig_*,\theta_{\sig_*-\tau}\w)$ for each $k\in\N$ such that
$$\phi_k=u_{\al_{n_k}}(\sig_k,\sig_*,\theta_{-\tau}\w,\tilde{\phi}_k).$$
Since $\tilde{\phi}_{k}\in\cup_{k\in\N^+}\cA_{\al_{n_k}}(\sig_*,\theta_{\sig_*-\tau}\w)$,
which is precompact in $\ell^2$ (by Theorem \ref{th3.4A} and Lemma \ref{le3.10}),
we can also assume that
\be\label{5.12A}\tilde{\phi}_k\ra\phi_*,\hs
\mb{as }n\ra\8.\ee
Recalling the condition (4) of a continuous NRDS in Section \ref{s2}
(or \cite[Theorem 2.5]{WJ24}) and the uniform continuity in Lemma \ref{le3.7},
we can infer from $\al_{n_k}\ra\al_0$, $\sig_k\ra\sig_0$ and $\tilde{\phi}_k\ra\phi_*$ that
\begin{align*}
\|\phi_k-u_{\al_0}(\sig_0,\sig_*,\theta_{-\tau}\w,\phi_*)\|\leqslant&
\|u_{\al_{n_k}}(\sig_k,\sig_*,\theta_{-\tau}\w,\tilde{\phi}_k)
-u_{\al_0}(\sig_k,\sig_*,\theta_{-\tau}\w,\tilde{\phi}_k)\|\\
&+\|u_{\al_{0}}(\sig_k,\sig_*,\theta_{-\tau}\w,\tilde{\phi}_k)
-u_{\al_0}(\sig_0,\sig_*,\theta_{-\tau}\w,\phi_*)\|\ra0,
\end{align*}
as $n\ra\8$.
We have proved the conclusion.
\eo

\bl\label{le5.7} Under the conditions of Theorem \ref{th5.5},
for each convergent sequence $\{\al_n\}_{n\in\N^+}$ in $\R$ with $\al_n\ra\al_0\in\R$,
and $\Psi\in\cT$, the family ($n\in\N^+$) of mappings
\be\label{5.22A}\vsig\mapsto\int_{H}\langle \tilde{F}(\vsig,u),\Psi'(u)\rangle
\mu^{\al_n}_{\vsig,\theta_{\vsig-\tau}\w}(\di u)\ee
is dominated by a Lebesgue integrable function over $[s,t]$.
\el
\bo By the compactness of $\ol{\cU[s,t]}$ by Lemma \ref{le5.6},
we know that $\|u\|$ is bounded on $\ol{\cU[s,t]}$.
Hence by \eqref{5.2a}, we have a positive constant $\sM$ independent of $s$ and $n$ such that
$$\left|\langle \tilde{F}(\vsig,u),\Psi'(u)\rangle\right|
\leqslant\sM(1+\|\psi_3(\vsig)\|_1+\|\psi_4(\vsig)\|_1)$$
and then
\begin{align}\label{5.23A}&\left|\int_{H}\langle \tilde{F}(\vsig,u),\Psi'(u)\rangle
\mu^{\al_n}_{\vsig,\theta_{\vsig-\tau}\w}(\di u)\right|=
\left|\int_{\cA_{\al_n}(\vsig,\theta_{\vsig-\tau}\w)}
\langle \tilde{F}(\vsig,u),\Psi'(u)\rangle
\mu^{\al_n}_{\vsig,\theta_{\vsig-\tau}\w}(\di u)\right|\notag\\
\leqslant&\sM(1+\|\psi_3(\vsig)\|_1+\|\psi_4(\vsig)\|_1),
\end{align}
which is Lebesgue integrable over $[s,t]$ by the assumption
$\psi_3,\psi_4\in L^1_{\rm loc}(\R,\ell^1)$.
We can immediately obtain the final result from \eqref{5.23A}.
\eo

\bl\label{le5.8} Under the conditions of Theorem \ref{th5.5},
for each convergent sequence $\{\al_n\}_{n\in\N^+}$ in $\R$ with $\al_n\ra\al_0\in\R$,
and $\Psi\in\cT$, the family ($n\in\N^+$) of mappings
\be\label{5.8A}\vsig\mapsto\int_{H}\Psi''(u)(u,u)
\mu^{\al_n}_{\vsig,\theta_{\vsig-\tau}\w}(\di u)\ee
is uniformly bounded over $[s,t]$.
\el
\bo Indeed by the condition (3) of the definition of $\cT$, we can easily see that
for all $u\in\ol{\cU[s,t]}\subset V$, $\Psi''(u)$ is bounded in $L(H\X H,\R)$.
Thus we obtain
$$\left|\Psi''(u)(u,u)\right|\leqslant\max_{u\in\ol{\cU[s,t]}}\|\Psi''(u)\|_{L(H\X H,\R)}\|u\|^2,$$
which is uniformly bounded by a positive constant $\sN$ independent of $s$ and $n$
on the compact set $\ol{\cU[s,t]}$. Therefore
$$\left|\int_{H}\Psi''(u)(u,u)\mu^{\al_n}_{\vsig,\theta_{\vsig-\tau}\w}(\di u)\right|
=\left|\int_{\cA_{\al_n}(\vsig,\theta_{\vsig-\tau}\w)}\Psi''(u)(u,u)
\mu^{\al_n}_{\vsig,\theta_{\vsig-\tau}\w}(\di u)\right|
\leqslant\sN.$$
The proof is complete.
\eo

Now we are ready to show Theorem \ref{th5.5} as follows.

\noindent\textit{Proof of Theorem \ref{th5.5}.}
By Theorem \ref{th4.4}, we have known that for each $\tau\in\R$, $\bP$-a.s. $\w\in\W$
and each convergent sequence $\{\al_n\}_{n\in\N^+}$ in $\R$ with $\al_n\ra\al_0\in\R$,
there exists a subsequence $\{n_k\}_{k\in\N^+}$ of $\{n\}_{n\in\N^+}$ such that
the invariant sample measures $\mu^{\al_{n_k}}_{\tau,\w}$ weakly converges to an
invariant sample measures $\mu^{\al_0}_{\tau,\w}$ of $\vp_{\al_0}$.
Then by definition it is natural that as $k\ra\8$,
$$\int_{H}\Psi(u)\mu^{\al_{n_k}}_{t,\theta_{t-\tau}\w}(\di u)\ra
\int_{H}\Psi(u)\mu^{\al_0}_{t,\theta_{t-\tau}\w}(\di u)\hs\mb{and}$$
$$\int_{H}\Psi(u)\mu^{\al_{n_k}}_{s,\theta_{s-\tau}\w}(\di u)
\ra\int_{H}\Psi(u)\mu^{\al_0}_{s,\theta_{s-\tau}\w}(\di u),$$
for the left sides of \eqref{5.11B} and \eqref{5.12B}.

For the right sides, the weak convergence of $\mu^{\al_{n_k}}_{\tau,\w}$
to $\mu^{\al_0}_{\tau,\w}$ ensures the convergence
$$\int_{H}\langle \tilde{F}(\vsig,u),\Psi'(u)\rangle
\mu^{\al_{n_k}}_{\vsig,\theta_{\vsig-\tau}\w}(\di u)\ra
\int_{H}\langle \tilde{F}(\vsig,u),\Psi'(u)\rangle
\mu^{\al_0}_{\vsig,\theta_{\vsig-\tau}\w}(\di u).$$
Combining the consequence in Lemma \ref{le5.7}, we can deduce from
the (Lebesgue) Dominated Convergence Theorem the convergence
$$\int_s^t\int_H\langle\tilde{F}(\vsig,u),\Psi'(u)\rangle
\mu^{\al_{n_k}}_{\vsig,\theta_{\vsig-\tau}\w}(\di u)\di\vsig
\ra\int_s^t\int_H\langle\tilde{F}(\vsig,u),\Psi'(u)\rangle
\mu^{\al_0}_{\vsig,\theta_{\vsig-\tau}\w}(\di u)\di\vsig.$$
The convergence
$$\frac{\al_{n_k}^2}{2}\int_s^t\int_H\Psi''(u)(u,u)
\mu^{\al_{n_k}}_{\vsig,\theta_{\vsig-\tau}\w}(\di u)\di\vsig
\ra\frac{\al_0^2}{2}\int_s^t\int_H\Psi''(u)(u,u)
\mu^{\al_0}_{\vsig,\theta_{\vsig-\tau}\w}(\di u)\di\vsig$$
is a similar consequence by Lemma \ref{le5.8}.

Now each term in \eqref{5.11B} except the term including stochastic integral converges to
the corresponding term in \eqref{5.12B}.
Moreover, the equations \eqref{5.11B} and \eqref{5.12B} are both obtained by Theorem \ref{th5.4}.
One can trivially infers that the remaining term including stochastic integral in \eqref{5.11B}
also converges as follows,
$$\al_{n_k}\int_s^t\int_H(u,\Psi'(u))\mu^{\al_{n_k}}_{\vsig,\theta_{\vsig-\tau}\w}(\di u)
\di\tilde{W}(\vsig)\ra\al_0\int_s^t\int_H(u,\Psi'(u))
\mu^{\al_0}_{\vsig,\theta_{\vsig-\tau}\w}(\di u)
\di\tilde{W}(\vsig),$$
as $k\ra\8$. This proves the termwise convergence now.\qed

\br As a special case of Theorem \ref{th5.5}, if we let $\al_0=0$,
then the equation \eqref{2.4} is a deterministic nonautonomous lattice equation.
Thus the convergence presented in Theorem \ref{th5.5} is actually from a stochastic Liouville type
equation to a deterministic Liouville type equation.
\er
\section{Convergence of sample statistical solutions}\label{s6}

In \cite{WZL24} the authors developed statistical solutions for deterministic equations
into sample statistical solutions for stochastic differential equations, which is defined
as follows in our framework.

\bd\label{de5.10} A family of Borel probability measures
$\{\mu^{\al}_{\tau,\w}\}_{(\tau,\w)\in\R\X\W}$
in $H$ is called a \textbf{sample statistical solution}
(in $H$) of \eqref{5.1}
if the following conditions are satisfied: for each $\tau\in\R$ and $\bP$-a.s. $\w\in\W$,
\benu\item[(1)] The function $s\mapsto\int_H\Psi(u)\mu^{\al}_{s,\theta_{s-\tau}\w}(\di u)$
is continuous for every $\Psi\in\cC_{\rm b}(H)$;
\item[(2)] For almost all $s\in\R$, the functions
$$u\mapsto\langle\tilde{F}(s,u),\phi\rangle,\hs u\mapsto(u,\phi)
\hs\mb{and}\hs u\mapsto\Phi(u,u)$$
are $\mu^{\al}_{s,\theta_{s-\tau}\w}$-integrable for every $\phi\in V$ and bilinear mapping
$\Phi\in L(H\X H,\R)$;
moreover, the mappings
$$s\mapsto\int_H\langle\tilde{F}(s,u),\phi\rangle\mu^{\al}_{s,\theta_{s-\tau}\w}(\di u)
\hs\mb{and}\hs s\mapsto\int_H\Phi(u,u)\mu^{\al}_{s,\theta_{s-\tau}\w}(\di u)$$
belongs to $L_{\rm loc}^1(\R)$, and the mapping
$$s\mapsto\int_H\(u,\phi\)\mu^{\al}_{s,\theta_{s-\tau}\w}(\di u)$$
belongs to $L^2_{\rm loc}(\R)$;
\item[(3)] For each $\Psi\in\cT$, the stochastic Liouville type equation
\eqref{1.4} holds.
\eenu
\ed

Then we can prove that the invariant sample measures obtained above are actually
a sample statistical solution for \eqref{5.1}.

\bt\label{th5.11} Suppose that the assumptions \textbf{(F1)}, \textbf{(F2)}, \textbf{(F3)},
\textbf{(F4)} and \textbf{(F5)} hold for $f$.
Then the family of Borel probability measures
$\{\mu_{\tau,\w}\}_{(\tau,\w)\in\R\X\W}$ obtained
in Theorem \ref{th4.2} is a sample statistical solution
of the equation \eqref{5.1}.
\et
\bo We check the conditions of Definition \ref{de5.10} one by one.
For the condition (1), we let $\Psi\in\cC_b(H)$.
Then by \eqref{4.15} and \eqref{4.17}, for $t,\,s,\,s_*\in\R$ with $t,\,s>s_*$
and $s,\,s_*$ fixed, we have
\begin{align}&\left|\int_{H}\Psi(u)\mu^{\al}_{t,\theta_{t-\tau}\w}(\di u)
-\int_H\Psi(u)\mu^{\al}_{s,\theta_{s-\tau}\w}(\di u)\right|\notag\\
\leqslant&\int_{\cA(s_*,\theta_{s_*-\tau}\w)}
\left|\Psi(u_\al(t,s_*,\theta_{-\tau}\w,u))
-\Psi(u_\al(s,s_*,\theta_{-\tau}\w,u))\right|
\mu^{\al}_{s_*,\theta_{s_*-\tau}\w}(\di u).\label{3.70z}
\end{align}
Since $\cA_\al(s_*,\theta_{s_*-\tau}\w)$ is compact,
by the continuity stated in (4) of the definition of NRDS's, we know the continuity of
$s\mapsto u_\al(s,s_*,\theta_{-\tau}\w,\phi)$ is uniform
for all $\phi\in\cA_\al(s_*,\theta_{s_*-\tau}\w)$.
Hence for every $\de>0$, there exists $\ol{\de}>0$
independent of $\phi\in\cA_\al(s_*,\theta_{s_*-\tau}\w)$, such that when $|t-s|<\ol{\de}$,
\be\label{3.71z}\|u_\al(t,s_*,\theta_{-\tau}\w,\phi)-u_\al(s,s_*,\theta_{-\tau}\w,\phi)\|<\de,\ee
for all $\phi\in\cA_\al(s_*,\theta_{s_*-\tau}\w)$.
Furthermore, for all $\phi\in\cA_\al(s_*,\theta_{s_*-\tau}\w)$ and $\vsig\in[s-1,s+1]$,
we know that
$$u_\al(\vsig,s_*,\theta_{-\tau}\w,\phi)\in\cU_\al[s-1,s+1],$$
which is compact by Lemma \ref{le5.1}.
Also, $\Psi\in\cC_{\rm b}(H)$ implies that $\Psi$ is uniformly continuous on $\cU_\al[s-1,s+1]$.
This means that, for each $\ve>0$, there is $\de>0$ such that when $u,\,v\in\cU_\al[s,t]$
and $\|u-v\|<\de$,
\be\label{3.72z}|\Psi(u)-\Psi(v)|<\ve.\ee

As a result, by \eqref{3.70z}, \eqref{3.71z} and \eqref{3.72z}, it can be seen that
for each $\ve>0$, there exists $\de>0$ and $\ol{\de}\in(0,1)$, such that when $|t-s|<\ol{\de}$,
\eqref{3.71z} holds true and then
\be\label{5.29z}\eqref{3.70z}<\int_{\cA_\al(s_*,\theta_{s_*-\tau}\w)}\ve
\mu^{\al}_{s_*,\theta_{s_*-\tau}\w}(\di\phi)=\ve,\ee
which indicates the condition (1) of Definition \ref{de5.10}.

For (2) of Definition \ref{de5.10},
we can recall Lemma \ref{le5.2}.
The conclusions for $\langle\tilde{F}(s,u),\phi\rangle$ and $(u,\phi)$
can be obtained by defining $\Psi(u)\equiv\phi$.
The conclusions for $\Phi(u,u)$ can be obtained by
defining $\Psi(u)=\Phi(u,u)$.

The condition (3) of Definition \ref{de5.10} follows from Theorem \ref{th5.4}.
The theorem is proved here.
\eo

According to the convergence consequences of invariant sample measures and
stochastic Liouville type equations,
we can incidentally obtain the convergence of sample statistical solutions.

\bt Suppose that the assumptions \textbf{(F1)}, \textbf{(F2)}, \textbf{(F3)},
\textbf{(F4)} and \textbf{(F5)} hold for $f$.
Let $\vp_\al$ be the NRDS defined as \ref{2.10} generated by \eqref{2.4}
and $\{\mu^{\al}_{\tau,\w}\}_{(\tau,\w)\in\R\X\W}$ be the sample statistical
solution of $\vp_\al$ given in Theorem \ref{th5.11}.
Then for each $\tau\in\R$, $\bP$-a.s. $\w\in\W$ and
each sequence $\{\al_n\}_{n\in\N^+}$ in $\R$
with $\al_n\ra\al_0$ as $n\ra\8$,
there exists a subsequence $\{\al_{n_k}\}_{k\in\N^+}$ of $\{\al_n\}_{n\in\N^+}$
depending on $\tau$ and $\w$,
and a sample statistical solution $\mu_{\tau,\w}^{\al_0}$ of $\vp_{\al_0}$ such that
$\mu_{\tau,\w}^{\al_{n_k}} $ converges to $\mu_{\tau,\w}^{\al_0}$ as $k\ra+\8$
in the following sense:
\benu\item[(1)]As a family of invariant sample measures, $\mu_{t,\theta_{t-\tau}\w}^{\al_{n_k}}$
converges weakly to $\mu_{t,\theta_{t-\tau}\w}^{\al_0}$ for every $t\in\R$;
\item[(2)]As $k\ra\8$,
\be\label{5.30z}\int_H\Psi(u)\mu^{\al_{n_k}}_{\vsig,\theta_{\vsig-\tau}\w}(\di u)\ra
\int_H\Psi(u)\mu^{\al_0}_{\vsig,\theta_{\vsig-\tau}\w}(\di u)\ee
in $\cC([s,t])$ with maximum norm for each interval $[s,t]\in\R$ and $\Psi\in\cC_{\rm b}(H)$;
\item[(3)]For all $\phi\in V$ and $\Phi\in L(H\X H,\R)$, as $k\ra\8$,
$$\int_H\langle\tilde{F}(\vsig,u),\phi\rangle\mu^{\al_{n_k}}_{\vsig,\theta_{\vsig-\tau}\w}(\di u)\ra
\int_H\langle\tilde{F}(\vsig,u),\phi\rangle\mu^{\al_0}_{\vsig,\theta_{\vsig-\tau}\w}(\di u)$$
$$\mb{and}\hs\int_H\Phi(u,u)\mu^{\al_{n_k}}_{\vsig,\theta_{\vsig-\tau}\w}(\di u)\ra
\int_H\Phi(u,u)\mu^{\al_0}_{\vsig,\theta_{\vsig-\tau}\w}(\di u)$$
in $L^1([s,t])$ and
\be\int_H\(u,\phi\)\mu^{\al_{n_k}}_{\vsig,\theta_{\vsig-\tau}\w}(\di u)\ra
\int_H\(u,\phi\)\mu^{\al_0}_{\vsig,\theta_{\vsig-\tau}\w}(\di u)\label{5.31z}\ee
in $L^2([s,t])$ for each interval $[s,t]\in\R$;
\item[(4)]As $k\ra\8$, the stochastic Liouville type equation \eqref{5.11B}
converges to \eqref{5.12B} termwise.
\eenu
\et

\bo The conclusion (1) follows from Theorem \ref{th4.4}.
We next apply Arzel\`a-Ascoli theorem to show the conclusion (2).
Observe that in the procedure of proving \eqref{5.29z},
as long as the compact set $\cU_{\al}[s,t]$ is replaced by $\ol{\cU[s,t]}$
(defined in Lemma \ref{le5.6}), which is compact,
\eqref{5.29z} is valid uniformly for all $\al=\al_{n}$, $n\in\N$.
This is indeed the equicontinuity of the sequence in \eqref{5.30z}.
The uniform boundedness is also obvious, since $\Psi$ has a uniform
bound on $\cU[s,t]$.
The reason why the limit integral is with respect to $\mu^{\al_0}_{\vsig,\theta_{\vsig-\tau}\w}$
is the weak convergence in (1).
Hence the convergence \eqref{5.30z} holds in $\cC[s,t]$
(perhaps for a subsequence).

The conclusion (3) follows from Lemmas \ref{le5.6} and \ref{le5.7},
where the convergence \eqref{5.31z} in $L^2([s,t])$ is a similar
result by the proof procedure of Lemma \ref{le5.8} with a slight modification.
The conclusion (4) has been proved in Theorem \ref{th5.5}.
\eo
\section*{Data Availability}

Data sharing is not applicable to this article as no new data were created
or analyzed in this study.

\section*{Acknowledgements}

Our work was supported by grants from the National Natural Science Foundation of China
(NNSFC Nos. 11801190 and 12101462).


\begin{thebibliography}{99}
	\small
	
	
	\bibitem{B21}B. Belean,
	Active contours driven by cellular neural networks
for image segmentation in biomedical applications,
	\emph{Stud. Inform. Control}, \textbf{30} (2021), 109-120.
	
\bibitem{CKV14}T. Caraballo, F. Morillas and J. Valero,
	On differential equations with delay in Banach spaces and attractors
for retarded lattice dynamical systems,
	\emph{Discrete Contin. Dyn. Syst.}, \textbf{34} (2014), 51-77.
	
	
\bibitem{CFZ23}P. Y. Chen, M. M. Freitas and X. P. Zhang,
	Random attractor, invariant measures, and ergodicity of
lattice $p$-Laplacian equations driven by superlinear noise,
	\emph{J. Geom. Anal.}, \textbf{33} (2023), Paper No. 98, 46 pp.

\bibitem{CY23}Z. Chen and D. D. Yang,
	Invariant measures and stochastic Liouville type theorem
for non-autonomous stochastic reaction-diffusion equations,
	\emph{J. Differential Equations}, \textbf{353} (2023), 225-267.

\bibitem{C02}I. Chueshov,
 \emph{Monotone Random Systems Theory and Applications}, Lect. Notes in Math., vol. 1779,
 Springer-Verlag, Berlin, 2002.
	
\bibitem{CLL18}H. Y. Cui, J. A. Langa and Y. R. Li,
  Measurability of random attractors for quasi strong-to-weak continuous random dynamical systems,
  \emph{J. Dyn. Differential Equaitons}, \textbf{30} (2018), 1873-1898.
	
\bibitem{CLY15}H. Y. Cui, Y. R. Li and J. Y. Yin,
 Existence and upper semicontinuity of bi-spatial pullback attractors for smoothing
 cocycles,
 \emph{Nonlinear Anal.}, \textbf{128} (2015), 303-324.
	
\bibitem{D77}K. Deimling,
	\emph{Ordinary Differential Equations in Banach Spaces}, Lect. Notes in Math., vol. 596,
	Springer-Verlag, Berlin, Heidelberg, New York, 1977.
	
\bibitem{FMRT01}C. Foias, O. Manley, R. Rosa and R. Temam,
 \emph{Navier-Stokes Equations and Turbulence},
 Cambridge University Press, Cambrige, 2001.	

\bibitem{GK16}A. H. Gu and P. E. Kloeden,
	Asymptotic behavior of a nonautonomous $p$-Laplacian lattice system,
	\emph{Int. J. Bifur. Chaos}, \textbf{26} (2016), Paper No. 1650174, 9 pp.
	
	
	\bibitem{GL17}A. H. Gu and Y. R. Li,
	Dynamic behavior of stochastic $p$-Laplacian-type lattice equations,
	\emph{Stoch. Dyn.}, \textbf{17} (2017), Paper No. 1750040, 19 pp.
	
	\bibitem{HSZ11}X. Y. Han, W. X. Shen and S. F. Zhou,
	Random attractors for stochastic lattice dynamical systems in weighted spaces,
	\emph{J. Differential Equations}, \textbf{250} (2011), 1235-1266.
	
	\bibitem{HLW21}Y. He, C. Q. Li and J. T. Wang,
	Invariant measures and statistical solutions
for the nonautonomous discrete modified Swift-Hohenberg equation,
	\emph{Bull. Malays. Math. Sci. Soc.}, \textbf{44} (2021), 3819-3837.

\bibitem{K05}H.-H. Kuo,
 \emph{Introduction to stochastic integration},
 Springer Science+Business Media, Inc., 2005.

\bibitem{LLW22}C. C. Li, C. Q. Li and J. T. Wang,
 Statistical solution and Liouville type theorem
 for coupled Schr\"odinger-Boussinesq equations on infinite lattices,
 \emph{Discrete Contin. Dyn. Syst. Ser. B}, \textbf{27} (2022), 6173-6196.
	
\bibitem{LLW23}C. C. Li, C. Q. Li and J. T. Wang,
 Statistical solution and Liouville type theorem for nonautonomous discrete Selkov model,
 \emph{Dyn. Syst.}, \textbf{38} (2023), 140-162.
	
\bibitem{LLF24}F. Z. Li, D. S Li and M. M. Freitas,
 Limiting dynamics for stochastic delay $p$-Laplacian equation on unbounded thin domains,
 \emph{Banach J. Math. Anal.}, \textbf{18} (13) (2024), 41 pp.

\bibitem{LLW19} F. Z. Li, Y. R. Li and R. H. Wang,
 Strong convergence of bi-spatial random attractors for parabolic equations
 on thin domains with rough noise,
 \emph{Topol. Methods Nonlinear Anal.}, \textbf{53} (2) (2019), 659-682.
	
\bibitem{LCL14}Y. R. Li, H. Y. Cui and J. Li,
 Upper semi-continuity and regularity of random attractors
 on $p$-times integrable spaces and applications,
 \emph{Nonlinear Anal.}, \textbf{109} (2014), 33-44.
		
	\bibitem{LGL15}Y. R. Li, A. H. Gu and J. Li,
	Existence and continuity of bi-spatial random attractors and application
to stochastic semilinear Laplacian equations,
	\emph{J. Differential Equaitons}, \textbf{258} (2015), 504-534.

	
	\bibitem{LFL14}P. M. Lima, N. J. Ford and P. M. Lumb,
	Computational methods for a mathematical model of propagation
of nerve impulses in myelinated axons,
	\emph{Appl. Numer. Math.}, \textbf{85} (2014), 38-53.
	
	
	\bibitem{N20}E. Nozawa,
	Coupled map lattice for the spiral pattern formation in astronomical objects,
	\emph{Phys. D}, \textbf{405} (2020), Paper No. 132377, 10 pp.
	
	\bibitem{SL22}T. F. Sequeria and P. M. Lima,
	Numerical simulations of one- and two-dimensional stochastic neural field equations with delay,
	\emph{J. Comput. Neurosci.}, \textbf{50} (2022), 299-311.
	
%
	\bibitem{SLW22}L. Song, Y. R. Li and F. L. Wang,
	Controller and asymptotic autonomy of random attractors
for stochastic $p$-Laplace lattice equations,
	\emph{Evol. Equ. Control Theory}, \textbf{11} (2022), 2033-2054.
	
	
%
	
	
\bibitem{W12}B. X. Wang,
 Sufficient and necessary criteria for existence of pullback attractors
 for non-compact random dynamical systems,
 \emph{J. Differential Equations}, \textbf{253} (2012), 1544-1583.
	
	\bibitem{W14}B. X. Wang,
	Existence and upper semicontinuity of attractors for stochastic equations
with deterministic non-autonomous terms,
	\emph{Stoch. Dyn.}, \textbf{14} (2014), Paper No. 1450009, 31 pp.
	
\bibitem{WJ24}J. T. Wang and W. H. Jin,
	Regularity of pullback random attractors and invariant sample measures
for nonautonomous stochastic $p$-Laplacian lattice system,
	\emph{Discrete Contin. Dyn. Syst. Ser. B}, \textbf{29} (3) (2024), 1344-1379.
	
\bibitem{WLYJ21}J. T. Wang, C. Q. Li, L. Yang and M. Jia,
Upper semi-continuity of random attractors and existence of invariant measures
for nonlocal stochastic Swift-Hohenberg equation with multiplicative noise,
	\emph{J. Math. Phys.}, \textbf{62} (2021), Paper No. 111507, 31 pp.
	
	
%
%
%
	
\bibitem{WZ23}J. T. Wang and X. Q. Zhang,
 Invariant sample measures and random Liouville type theorem
 for a nonautonomous stochastic $p$-Laplacian equation,
 \emph{Discrete Contin. Dyn. Syst. Ser. B}, \textbf{28} (2023), 2803-2827.
	
\bibitem{WZL23}J. T. Wang, X. Q. Zhang and C. Q. Li,
 Global martingale and pathwise solutions and infinite regularity of invariant measures
 for a stochastic modified Swift-Hohenberg equation,
 \emph{Nonlinearity}, \textbf{36} (2023), 2655-2707.
	
\bibitem{WZZ21}J. T. Wang, X. Q. Zhang and C. D. Zhao,
 Statistical solutions for a nonautonomous modified Swift-Hohenberg equation,
 \emph{Math. Methods Appl. Sci.}, \textbf{44} (2021), 14502-14516.
	
\bibitem{WZC20}J. T. Wang, C. D. Zhao and T. Caraballo,
 Invariant measures for the 3D globally modified Navier-Stokes equations
 with unbounded variable delays,
 \emph{Commun. Nonlinear Sci. Numer. Simul.}, \textbf{91} (2020), Paper No. 105459, 14 pp.

\bibitem{WZL24}J. T. Wang, D. D. Zhu and C. Q. Li,
 Invariant sample measures and sample statistical solutions for
 nonautonomous stochastic lattice Cahn-Hilliard equation with nonlinear noise,
 \emph{arXiv}, arXiv: 2404.14798, 48 pp.	

\bibitem{WCT23}R. H. Wang, T. Caraballo and N. Tuan,
 Asymptotic stability of evolution systems of probability measures for nonautonomous stochastic systems: Theoretical results and applications,
 \emph{Proc. Amer. Math. Soc.}, \textbf{151} (2023), 2449-2458.
	
\bibitem{WW20}R. H. Wang and B. X. Wang,
 Random dynamics of $p$-Laplacian lattice systems driven by infinite-dimensional nonlinear noise,
 \emph{Stochastic Process Appl.}, \textbf{130} (2020), 7431-7462.
	
\bibitem{WWL22}X. J. Wang, J. T. Wang and C. Q. Li,
 Invariant measures and statistical solutions
 for a nonautonomous nonlocal Swift-Hohenberg equation,
 \emph{Dyn. Syst.}, \textbf{37} (2022), 136-158.

\bibitem{WWHK}Y. J. Wang, Y. Wang, X. Y. Han and P. E. Kloeden,
 A two-dimensional stochastic fractional non-local diffusion lattice model with delays,
 \emph{Stoch. Dyn.}, \textbf{22} (2023), Paper No. 2240032.
	
\bibitem{XZL20}C. B. Xiu, R. X. Zhou and Y. X. Liu,
 New chaotic memristive cellular neural network and
 its application in secure communication system,
 \emph{Chaos Solitons Fractals}, \textbf{141} (2020), Paper No. 110316, 15 pp.

\bibitem{XC22}J. H. Xu and T. Caraballo,
 Long Time Behavior of Stochastic Nonlocal Partial Differential Equations and Wong-Zakai Approximations,
 \emph{SIAM J. Math. Anal.}, \textbf{54} (2022), 2792-2844.

\bibitem{YCC23}D. D. Yang, T. Caraballo and Z. Chen,
 The periodic and limiting behaviors of invariant measures for 3D globally modified
 Navier-Stokes equations,
 \emph{J. Dynam. Differential Equations}, (2023), 21 pp.
 https://doi.org/10.1007/s10884-023-10260-8
	
\bibitem{YZ24}C. W. Yu and Y. Zeng,
 Local second order Sobolev regularity for $p$-Laplacian equation in semi-simple Lie group,
 \emph{Mathematics}, \textbf{12} (4) (2024), 14 pp.
	
\bibitem{ZZ23}X. H. Zhang and X. P. Zhang,
 Upper semi-continuity of non-autonomous fractional stochastic $p$-Laplacian equation
 driven by additive noise on $\R^n$,
 \emph{Discrete Contin. Dyn. Syst. Ser. B}, \textbf{28} (1) (2023), 385-407.
	
\bibitem{ZWC22}C. D. Zhao, J. T. Wang and T. Caraballo,
 Invariant sample measures and random Liouville type theorem
 for the two-dimensional stochastic Navier-Stokes equations,
 \emph{J. Differential Equaitons}, \textbf{317} (2022), 474-494.
	
\end{thebibliography}
\end{document}